\documentclass[11pt]{article}
\title{Universality results for largest eigenvalues of some sample covariance matrix ensembles}
\author{Sandrine P\'ech\'e \footnote{Permanent address: Institut Fourier BP 74, 100 Rue des
    maths, 38402 Saint Martin d'Heres, France; E-mail: Sandrine.Peche@ujf-grenoble.fr Temporary address: Department of Mathematics,
University of California at Davis, 
One Shields Ave., Davis, CA 95616, USA.}}
\usepackage{epsfig}
\usepackage{amsfonts}
\usepackage{amssymb}
\usepackage{amsthm}
\usepackage{amstext}
\usepackage{amscd}
\usepackage{amsmath}
\usepackage{setspace}
\usepackage{delarray}
\begin{document}
\newtheorem{theo}{Theorem}[section]
\newtheorem{prop}{Proposition}[section]
\newtheorem{lemme}{Lemma}[section]
\newtheorem{conjecture}{Conjecture}[section]
\newtheorem{coro}{Corollary}[section]
\newtheorem{definition}{Definition}[section]
\newtheorem{fact}{Fact}[section]
\newtheorem{hyp}{Assumption}[section]
\theoremstyle{remark}
\newtheorem{rem}{Remark}
\newtheorem{remark}{Remark}[section]
\newtheorem{Remark}{Remark}[section]
\newtheorem{Notationnal remark}{Remark}[section]
\newcommand{\bremnot}{\begin{Notationnal remark}}
\newcommand{\eremnot}{\end{Notationnal remark}}
\newcommand{\brem}{\begin{remark}}
\newcommand{\erem}{\end{remark}}
\newcommand{\bconj}{\begin{conjecture}}
\newcommand{\econj}{\end{conjecture}}
\newcommand{\bdefi}{\begin{definition}}
\newcommand{\edefi}{\end{definition}}
\newcommand{\bt}{\begin{theo}}
\newcommand{\bfa}{\begin{fact}}
\newcommand{\efa}{\end{fact}}
\newcommand{\Si}{\Sigma}
\newcommand{\becor}{\begin{coro}}
\newcommand{\ecor}{\end{coro}}
\newcommand{\et}{\end{theo}}
\newcommand{\bp}{\begin{prop}}
\newcommand{\ep}{\end{prop}}
\newcommand{\bl}{\begin{lemme}}
\newcommand{\el}{\end{lemme}}
\newcommand{\be}{\begin{equation}}
\newcommand{\ee}{\end{equation}}
\newcolumntype{L}{>{$}l<{$}}
\newenvironment{Cases}{\begin{array}\{{lL.}}{\end{array}}
\maketitle

\begin{abstract}
For sample covariance matrices with iid entries with sub-Gaussian 
tails, when both the number of samples and the number of variables 
become large and the ratio approaches to one, it is a well-known result 
of A. Soshnikov that the limiting distribution of the largest eigenvalue is 
same as the of Gaussian samples. In this paper, we extend this result 
to two cases. The first case is when the ratio approaches to an 
arbitrary finite value. The second case is when the ratio becomes infinity or arbitrarily small.\end{abstract}

%\begin{keyword}[class=AMS]
%\kwd[Primary ]{}
%\kwd{}
%\kwd[; secondary ]{}
%\end{keyword}

%\begin{keyword}
%\kwd{}
%\kwd{}
%\end{keyword}

\section{Introduction}
The scope of this paper is to study the limiting behavior of the largest eigenvalues of real
and complex sample covariance matrices with independent identically distributed (i.i.d.), but non necessarily Gaussian, entries.
Consider a sample of size $p$ of i.i.d. $N \times 1$ random vectors
$\vec{y}_1,\dots, \vec{y}_p$. We further assume that the sample vectors $\vec{y}_k$ have mean zero and \emph{covariance} $\Si=Id$.
We use $X=[\vec{y}_1,\cdots, \vec{y}_p]$ to denote the $N\times p$ data matrix and $M_N=\frac1{N} XX^*$ to denote
the \emph{sample covariance matrix}.
Random sample covariance matrices have been first studied in mathematical statistics
(\cite{Bronk}, \cite{James}, \cite{Hotelling}). A huge literature deals with the case where $p\to \infty$, $N$ being fixed, which is now quite well understood. Contrary to the traditional assumptions, it is now of current interest to study the case where $N$ is of the same order as $p$, due to the large amount of data available.
In particular, the limiting behavior of the largest eigenvalues is important for testing hypotheses on the covariance matrix $\Sigma.$ Here we focus on the simple case, $H_o: \Sigma=Id$ versus $H_a: \Sigma \not=Id,$ and study the asymptotic distribution of extreme eigenvalues under the $H_o$.
The study of extreme eigenvalues is also of interest in principal component analysis. We refer the reader to \cite{Johnstone}
and \cite{NEKlarginterest} for a review of statistical applications.
Other examples of applications include genetics \cite{Paterson},
mathematical finance \cite{PlerousGRAGS}, \cite{LalouxCPB},
\cite{MalevergneS}, wireless communication \cite{Telatar}, physics
of mixture \cite{SearC}, and statistical
learning \cite{HoyleR}.
We point out that the spectral properties of $M_N$ readily translate to the companion matrix $W_N=\frac{1}{N}X^{*}X.$ Indeed, $W_N$ is a $p\times p$ matrix, of rank $N$, with the same non-zero eigenvalues as $M_N.$ Thus, it is enough to study the spectral properties of $M_N$ to give a complete picture of the spectrum of such sample covariance matrices.

\subsection{Model and results}
We consider both real and complex random sample covariance matrices $$M_N=\frac{1}{N}XX^*,$$
where $X$ is a $N\times p,$ $p=p(N)\geq N$, random matrix satisfying certain ``moment conditions''. In the whole paper, we set $\gamma_N=\frac{p}{N}.$ 
We assume that the entries $X_{ij}, 1\leq i\leq N, 1\leq j\leq p,$ of the sequence of random matrices $X=X_N$ are non-necessarily Gaussian random variables satisfying the following conditions. First, in the complex case, \\
(i)$ \{ \Re eX_{i,j} , \, \Im mX_{i,j}: \:  1 \leq i  \leq N,\: 1\leq j\leq p\}$ are real independent random
variables,\\
(ii) all these real variables have symmetric laws (thus,
$\mathbb E[X_{i,j}^{2k+1}] =0$ for all $k \in \mathbb N ^*$),\\
(iii) $ \forall \, i ,j, \, \mathbb E[(\Re e X_{i,j})^{2}]=\mathbb E[(\Im mX_{i,j})^{2}]=
\frac{\sigma^2}{2}$,\\
(iv) all their other moments are assumed to be sub-Gaussian i.e. there exists a
constant $\tau >0$ such that uniformly in $i,j$ and $k$,
$$\mathbb E[|X_{i,j}|^{2k}] \leq {(\tau  \, k )}^k.$$
In the real setting, $X=(X_{i,j})_{ 1 \leq i , j \leq N}$ is a random matrix
such that\\
(i') the $\{ X_{i,j}, \: 1 \leq i  \leq N\,, 1\leq j\leq p\}$ are independent random variables,\\
(ii') the laws of the $X_{i,j}$ are symmetric (in particular,
$\mathbb E[X_{i,j}^{2k+1}] =0$),\\
(iii') for all $i ,j$, $\mathbb E [X_{i,j} ^2]=\sigma^2$, \\
(iv') all the other moments of the $X_{i,j}$ grow not faster than the Gaussian ones. This
means that there is a constant $\tau >0$ such that, uniformly in $i, \, j$
and $k$, $\, \mathbb E[X_{i,j}^{2k}] \leq {(\tau  \, k )}^k$.\\
When the
entries of $X$ are further assumed to be Gaussian, we will denote by $X_G$ the
corresponding model. In this case and in the complex setting,  $M_N^G$ is of the
so-called Laguerre Unitary Ensemble (LUE), which is also called the complex Wishart ensemble. In the real setting, $M_N^G$ is of the so-called Laguerre Orthogonal Ensemble (LOE) or real Wishart ensemble.

\paragraph{}The scope of this paper is to describe the large-$N$-limiting distribution of the $K$ largest eigenvalues induced by any such ensemble, for any fixed integer $K$ independent of $N.$ Two regimes are investigated in this paper. In the first part, we assume that there exists some constant $\gamma\geq 1$ such that $\lim_{N \to \infty} \gamma_N =\gamma.$ 
In the second part, we consider the case where $\gamma_N\to \infty$ as $N \to \infty.$\\ 
Before stating our results, we recall some known results about sample covariance matrices. We first focus on the case where $\lim_{N \to \infty}\gamma_N=\gamma<\infty.$
Let $\lambda_1\geq \lambda_2\geq \cdots \geq \lambda_N$ be the ordered eigenvalues induced
by any ensemble of the above type. The first fundamental result for the limiting spectral behavior
 of such random matrix ensembles has been obtained by Marchenko and Pastur in \cite{MP} (in a much more general context than here).
It is in particular proved therein that the spectral measure $\mu_{N}=\frac{1}{N}\sum_{i=1}^N\delta_{\lambda_i}$ a.s. converges as $N $ goes to infinity. Set $u_{\pm}^c=\sigma^2(1\pm\sqrt \gamma)^2.$ Then one has that 
\be \lim_{N \to \infty }\mu_N =\rho_{MP}\text{ a.s.}, \text{ where }\frac{d\rho_{MP}(x)}{dx}=\frac{\sqrt{(u_+^c-x)(x-u_-^c)}}{2\pi x \sigma
^2}1_{[u_-^c, u_+^c]}(x).\label{MarP}\ee  The limiting probability distribution $\rho_{MP}$ is the so-called Marchenko-Pastur distribution.\\ 
The above result gives no insight
about the behavior of the largest eigenvalues. The first study of the asymptotic behavior of the
largest eigenvalue goes back to S. Geman \cite{Geman}. It was later refined in \cite{BaiYinKrish} and
\cite{Silversteinwklim}. In particular, it is well known that $\lim_{N \to \infty} \lambda_1=u_+^c$
a.s. if the entries of the random matrix $X$ admit moments up to order $4$. 
Significant results about fluctuations of the largest eigenvalues around $u_+^c$ are much more
recent and are essentially established for Wishart ensembles only. In particular, the limiting distribution of the largest eigenvalue has been obtained by K. Johansson \cite{JohanssonWIs} for complex Wishart matrices and I. Johnstone \cite{Johnstone} for real Wishart matrices. A. Soshnikov \cite{SosWish} has derived for both ensembles the limiting distribution of the $K$ largest eigenvalues, for any fixed integer $K$.
Before recalling their results, we need a few definitions.
We denote by $\lambda_1^{G,\beta}\geq \lambda_2^{G,\beta}\geq \cdots \geq \lambda_N^{G,\beta}$ the
eigenvalues induced by the Wishart ensembles, with $\beta=2$ (resp. $\beta=1$) for the LUE (resp.
the LOE). We also define the limiting Tracy-Widom distribution for the largest eigenvalue. 
Let $Ai$ denote the standard Airy function and $q$ denote the solution of the Painlev\'e II differential equation $\frac{\partial ^2}{\partial x^2}q=xq(x)+2q^3(x),$ with boundary condition $q(x)\sim Ai(x)\text{ as }x \to +\infty.$
\bdefi The GUE (resp. GOE) Tracy-Widom distribution for the largest eigenvalue is defined by the cumulative distribution function $F_2(x)$ $=\exp{\{\int_x^{\infty}(x-t)q^2(t)dt\}}$ (resp. $F_1(x)=\!\exp{\{\int_x^{\infty}\frac{-q(t)}{2}+\frac{(x-t)}{2}q^2(t)dt\}}).$
\edefi %
The GUE (resp. GOE) Tracy-Widom distribution for the joint distribution of the $K$ largest eigenvalues (for any fixed integer $K$) has been similarly defined. We refer the reader to \cite{TWAi} and \cite{TWAi2} for a precise definition. \\
We then rescale the eigenvalues as follows: for $i=1, \ldots, N,$ we set %
\be \label{rescaling} \mu_i^{(G, \beta)}= \frac{\gamma_N^{1/6}}{(1+\sqrt{\gamma_N})^{4/3}}\frac{N^{2/3}}{\sigma^2} \left ( \lambda_i^{(G, \beta)}
-\sigma^2(1+\sqrt{\gamma_N})^2\right).\ee

 \bt\label{theo Johns}
\cite{JohanssonWIs} \cite{Johnstone} \cite{SosWish}. The joint distribution of the $K$ largest
eigenvalues of the LUE (resp. LOE) rescaled as in (\ref{rescaling}) converges, as $N \to \infty$,
to the joint distribution defined by the GUE (resp. GOE) Tracy-Widom law.
\et 
The proof of Theorem \ref{theo Johns} relies on the crucial fact that the joint eigenvalue density of the
Wishart ensembles can be exactly computed. 
Starting from numerical simulations, it was then conjectured, in \cite{Johnstone} e.g., that
Theorem \ref{theo Johns} actually holds for a class of random sample covariance matrices much wider than the Wishart ensembles. 
Such a universality result was later proved for some quite general ensembles by  A. Soshnikov \cite{SosWish}, yet
under some restriction on the sample size, as we now recall. \\
For any ensemble satisfying $(i)$ to $(i\nu)$ (resp. $(i')$ to $(i\nu')$), we set:
\be \label{rescalinggen} \mu_i= \frac{\gamma_N^{1/6}}{(1+\sqrt{\gamma_N})^{4/3}}\frac{N^{2/3}}{\sigma^2} \left ( \lambda_i-\sigma^2(1+\sqrt{\gamma_N})^2\right).\ee

\bt \cite{SosWish}\label{theo: unisos} Assume that $p-N=O(N^{1/3})$. The joint distribution of the $K$ rescaled largest eigenvalues $\mu_i, i\leq K$, induced by any ensemble satisfying $(i)$ to
$(i\nu)$ (resp. $(i')$ to $(i\nu')$) converges, as $N$ goes to infinity, to the joint distribution
defined by the GUE (resp. GOE) Tracy-Widom law. \et

In this paper, we prove that such a universality result holds for any value of the parameter $\gamma.$ This is the main result of this note.

\bt \label{theo: unigammafini} The joint distribution of the $K$ rescaled largest eigenvalues $\mu_i, i\leq K$, induced by any ensemble satisfying $(i)$ to $(i\nu)$ (resp. $(i')$ to $(i\nu')$) converges, as $N$ goes to infinity, to the joint distribution
defined by the GUE (resp. GOE) Tracy-Widom law. The results holds for any value of the parameter $\gamma\geq 1.$ \et
\brem Assumptions $(i\nu)$ and $(i\nu')$ can actually be relaxed. This relaxation is discussed in the second paragraph of Subsection \ref{Subsec: stat}. 
\erem
\paragraph{}Before giving secondary results, we give a few comments on the way we proceed to prove Theorem \ref{theo: unigammafini}. In Theorem \ref{theo: unisos}, the reason for the restriction on $p-N$ follows from the idea of the proof used therein. Basically,
when $\gamma=1$, the eigenvalues of a random sample covariance matrix roughly behave as the squares of those of a typical Wigner random matrix. This adequacy still works for the largest eigenvalues,
but fails if $\gamma$ is not close enough to one. Theorem \ref{theo: unisos} has been proved using universality
results established for classical Wigner random matrices.
Here, we revisit the problem of computing the asymptotics of $\mathbb{E}\left[\text{Tr} M_N^L\right]$ for some powers $L$ that may go to infinity, using combinatorial tools specifically well suited for the study of spectral functions of sample covariance matrices. It is well known that Dyck paths and Catalan numbers are associated to standard Wigner matrices (see \cite{BaiMethods}). Suitable combinatorial tools in the case of sample covariance matrices are the so-called \emph{Narayana numbers} and some particular Dyck paths. Using those, we can extend the universality result of \cite{SosWish} to any value of the ratio $\gamma.$ \\
The case where $\gamma\leq 1$ can also be considered thanks to the companion matrix $W_N$.
Let $\lambda'_i,1 \leq i \leq p,$ be the eigenvalues of $\frac{1}{p}X^*X$, ordered in decreasing order and let $\delta_N=\gamma_N^{-1}=N/p,$ so that $\delta_N \to 1/\gamma\leq 1$ as $N \to \infty.$ We set:
$$\mu'_i=\frac{{\delta_N}^{1/6}}{(1+\sqrt{{\delta_N}})^{4/3}}p^{2/3}\left (\frac{\lambda'_1}{{\sigma^2}}-(1+\sqrt{{\delta_N}})^2\right ), i=1, \ldots,p.$$  
\becor Under the assumptions $(i)$ to $(i\nu)$ (resp. $(i')$ to $(i\nu')$),
the joint distribution of $\left (\mu'_1,\ldots,\mu'_K\right )$ converges as $N \to \infty$ to the GUE (resp. GOE) Tracy-Widom joint distribution of the $K$ largest eigenvalues. 
\ecor

The machinery we develop to prove Theorem \ref{theo: unigammafini} can also be  used to consider the case where the size of the sample data increases in such a way that $p, N$ go to infinity and $ \frac{p}{N}\to \infty$. The large-$N$-limiting behavior of extreme eigenvalues of Wishart matrices for such a regime has been obtained by N. El Karoui \cite{NEK}. The particular interest of such a study for statistical applications (e.g. microarrays) is also explained in great detail therein.
%The first results for such a regime have been studied for Wishart ensembles by N. El Karoui \cite{NEK}. Therein the almost sure limit of the spectral measure of such random matrices is investigated, as well as the large-$N$-limiting behavior of extreme eigenvalues. The particular interest of such a study for statistical applications (e.g. microarrays) is also explained in great detail therein.
\bt \label{theoNEK} \cite{NEK} With the same rescaling as in (\ref{rescaling}), Theorem \ref{theo
Johns} actually holds in the case where $\lim_{N \to \infty}\gamma_N= \infty$. \et

Under the same assumptions $(i)$ to $(iv)$ (resp. $(i')$ to $(iv')$), we prove that universality
still holds in the regime $p/N \to \infty.$

\bt  \label{theo: unigammainfini}Theorem \ref{theo: unigammafini} also holds if
$\lim_{N \to \infty}\gamma_N= \infty$. \et
\subsection{Statistical implications of the result \label{Subsec: stat}}Testing homogeneity of a population has long been of interest in mathematical statistics, and it is often a preliminary step in discriminant analysis and cluster analysis.
Assuming high dimensionality, we consider the test of the null hypothesis $H_o:$ $\Sigma=Id$ vs. the alternative hypothesis $H_a: \Sigma\not=Id.$ The result stated here for the largest eigenvalue can be formulated as \be \label{test}\lim_{N \to \infty}P(\mu_1\leq x|H_o)=F_{2(1)}(x).\ee
The above theoretical result was essentially established for Gaussian samples only so far (cf. \cite{JohnstoneICM} for a review). Removing the Gaussianity assumption is actually fundamental for various statistical problems. 
Our result may be of use for instance in genetics (see \cite{Paterson} e.g.). Samples in genetic data are usually drawn from a distribution with compact support and the size of matrices encountered therein is typically large enough so that (\ref{test}) should be observed for some appropriate models. Some Gaussian (or other kinds) mixtures also fall into the class of distributions studied here. Such distributions occur for instance in finance in modeling some fat-tailed returns.
\paragraph{}Regarding especially the assumptions we make on the distribution of the entries, they may appear strong for other statistical purposes. The moment assumptions $(i\nu)$ and $(i\nu')$ can actually be relaxed, using truncation techniques. We can show that Theorem \ref{theo: unigammafini} holds under the assumption that 
$P (|X_{ij}|>x)\leq C(1+x)^{m_o}, \forall i,j,$ for some $m_o>36$ (see Remarks \ref{Rem: cutoff} and \ref{Rem: cutoffmaj}). 
We do not consider this case here, which would increase the technicalities of the paper. To illustrate this, a simulation is given below where the entries of $X$ have a Student's  $t-$distribution with $40$ degrees of freedom. These distributions may be interesting in statistical models due to their (relative) robustness with respect to outliers.\\
The symmetry assumption is probably more problematic and is again a technical assumption for the proof.
Indeed, it is expected that the lack of symmetry has no impact on the limiting distribution of largest eigenvalues (provided the distribution is centered). Yet, analytical tools to prove such a result are not established (see e.g. \cite{PecheSos} for recent progress).
\paragraph{}
The method we develop is also a first step towards considering samples with non-Identity covariance. Such results are of practical importance for understanding the behavior of Principal Component Analysis and dimension reduction in high dimensional setting. It is therefore important to consider covariance matrices with more complex structure. In particular (in progress), the moment approach developed here seems to be well suited in the case where the population covariance is a so-called ``spiked'' diagonal matrix. That is, $\Sigma= Id + D$, where the deformation $D$ is a finite rank diagonal matrix. This is important, since the test based on (\ref{test}) may not reject $H_o$ if the largest eigenvalue of $D$ is not large enough, because of a phase transition phenomenon described e.g. in \cite{BaikSilverstein}.

\paragraph{}A few simulations have been done to give, from a practical point of view, an idea of the rate of convergence of the distribution of the largest eigenvalue. We have generated real random matrices with i.i.d. entries with a $t-$distribution or a Gaussian mixture distribution. To fit the limiting Tracy-Widom distribution, we have rescaled the largest eigenvalue as follows:
\be \lambda_1^{N,p}:= \frac{N\lambda_1 -\sigma^2( \sqrt{N+a_1}+\sqrt{p+a_2})^2}{\sigma^2( \sqrt{N+a_1}+\sqrt{p+a_2} )\left( \frac{1}{\sqrt{N+a_1}}+\frac{1}{\sqrt{p+a_2}}\right )^{1/3}},\label{ressimul}\ee
for some adjustment parameters $a_1$ and $a_2$. We indeed have some freedom in the choice of these parameters, which can be any fixed real numbers. In the real Wishart case, the best parameters are known to be $a_1=a_2=-0.5.$ Determining the optimal parameters is important to improve convergence rates for (\ref{test}). They have been established in \cite{NEKCR} for complex Gaussian samples. Providing a general formula for these parameters is an issue that we cannot handle so far. Nevertheless, in view of our simulations, the optimal parameters may depend on the distribution of the entries and maybe also on the dimensions $N$ and $p.$ For instance, for the Gaussian mixture distribution, we found empirically that the best parameters are $a_1=a_2=0$ (choosing $a_1=a_2=-0.5$ gives results which are similar but not so satisfying). We also tested various $t-$distributions. For the $t-$distribution with $40$ (resp. $20$) degrees of freedom, we found that $a_1=a_2=-0.5$ (resp. $a_1=a_2=0$) were the optimal parameters.
Such a change can be understood, as the similarity between the $t-$ distribution and the Gaussian distribution increases with the number of degrees of freedom.\\
For the simulations given below, we have considered two distributions:\\
- the Gaussian mixture  $\frac{1}{2}\mathcal{N}(0,1)+\frac{1}{2}\mathcal{N}(0,3)$\\
- the Student's $t-$distribution with $40$ degrees of freedom.\\
We also considered various dimensions $N$, $N=10,25$, and $50$, as well as different sample size to dimension ratios, $\gamma_N=2,4,50,$ and $100$. \\
For each size and each distribution, we have generated $R=10000$ random matrices $X$ with i.i.d. entries. For each replication, we have rescaled the largest eigenvalue of $XX^*/N$ as in (\ref{ressimul}) with $a_1=a_2=-0.5$ (resp. $a_1=a_2=0$) for the Student's $t-$ (resp. Gaussian mixture) distribution.
We then derived the estimated cumulative probabilities for $\lambda_1^{N,p}$ obtained from the $10000$ replications. \\
For small sizes, $N=10$ and $p=20$ e.g, the proximity between the observed cumulative distribution of the rescaled largest eigenvalue and the Tracy-Widom distribution is reasonable essentially for upper quantiles for the $t-$distribution(95 $\%$). For the Gaussian mixture, it is reasonable for smaller quantiles.
As the size of the matrix increases, the proximity becomes acceptable for almost the whole range. We also note that the convergence is almost as good as for the Wishart ensemble (compare with Table 1 in \cite{Johnstone}) in both cases.  

\begin{table}\caption{Gaussian mixture, $\gamma_N=2,4$}
\begin{center}\scriptsize
\begin{tabular}{c|ccccccc} 
 $P_c$. & $F_1$ & $10\times 20$ & $25\times 50$ & $50 \times 100$ & $10\times 40$  & $25\times 100$  & $50 \times 200$\\
\hline
-3.896	&.01	&.0079	 &.0103  &.0097  &.0107  &.0114	& .0092\\
-3.516	&.025	&.0256	 &.0246  &.0270  &.0257  &.0267	& .0236\\
-3.180	&.05	&.0554	 &.0487  &.0512  &.0568  &.0516	& .0489 \\
-2.782	&.10	&.1172	 &.0973  &.0971  &.1099  &.1004	& .0994\\
-2.088	&.25	&.2912	 &.2502  &.2478  &.2724  &.2453	& .2436\\
-1.269	&.50	&.5354	 &.4935  &.4894  &.5154  &.4960	& .4922\\
-0.392  &.75	&.7550	 &.7401  &.7336  &.7465  &.7440	& .7389\\
0.450	&.90	&.8951	 &.8897  &.8879  &.8870  &.8892	& .8894\\
0.979	&.95	&.9417	 &.9396  &.9436  &.9368  &.9422	& .9425\\
1.454	&.975	&.9676	 &.9662  &.9718  &.9644  &.9718	& .9691\\
2.024	&.99	&.9855	 &.9857  &.9875  &.9824  &.9875	& .9858\\
\end{tabular}
\end{center}
{\small
The first column shows the percentiles of the $F_1$ Tracy-Widom distribution corresponding to the values in the second column. The next $6$ columns give the estimated cumulative probabilities for $\lambda_1^{N,p}$ obtained from $10000$ replications. The entries of the random matrices are i.i.d. with the Gaussian mixture distribution $1/2 \mathcal{N}(0,1)+1/2 \mathcal{N}(0,3)$ and $a_1=a_2=0$.
The \texttt{Matlab} functions \texttt{normrnd} and \texttt{unifrnd} were used to generate the Gaussian mixtures. The Tracy-Widom quantiles were computed thanks to the \textbf{p2Num} package, provided by C. Tracy.}
\end{table}
\begin{table}\caption{Student's $t-$distribution, $\gamma_N=2,4$}
\begin{center}\scriptsize
\begin{tabular}{c|ccccccc} 
 $P_c$. & $F_1$ & $10\times 20$ & $25\times 50$ & $50 \times 100$ & $10\times 40$  & $25\times 100$  & $50 \times 200$\\
\hline
-3.896	&.01	&.0024	 &.0047  &.0069 &.0041  &.0052  & .0078 \\
-3.516	&.025	&.0097	 &.0144  &.0188 &.0129  &.0162  & .0220\\
-3.180	&.05	&.0262	 &.0346  &.0392 &.0336  &.0350  & .0469\\
-2.782	&.10	&.0674	 &.0838  &.0857 &.0755  &.0842  & .0924\\
-2.088	&.25	&.2203	 &.2266  &.2375 &.2220  &.2317  & .2399\\
-1.269	&.50	&.4863	 &.4784  &.4950 &.4835  &.4890  & .4922\\
-0.392  &.75	&.7457	 &.7421  &.7487 &.7415  &.7433  & .7493 \\
0.450	&.90	&.8908	 &.8950  &.9028 &.8951  &.8888  & .9003\\
0.979	&.95	&.9440	 &.9463  &.9501 &.9490  &.9421  & .9507\\
1.454	&.975	&.9694	 &.9730  &.9730 &.9727  &.9697  & .9735\\
2.024	&.99	&.9872	 &.9887  &.9896 &.9878  &.9880  & .9895\\
\end{tabular}
\end{center}
{\small
The entries of the random matrices are i.i.d. with a $t-$distribution with $40$ degrees of freedom and $a_1=a_2=-0.5$. We have used the \texttt{Matlab} function \texttt{trnd} to generate the Student random variables.}
\end{table}

\begin{figure}
\begin{center}
\begin{tabular}{cp{.5cm}c}
\includegraphics[height=4cm,angle=0]{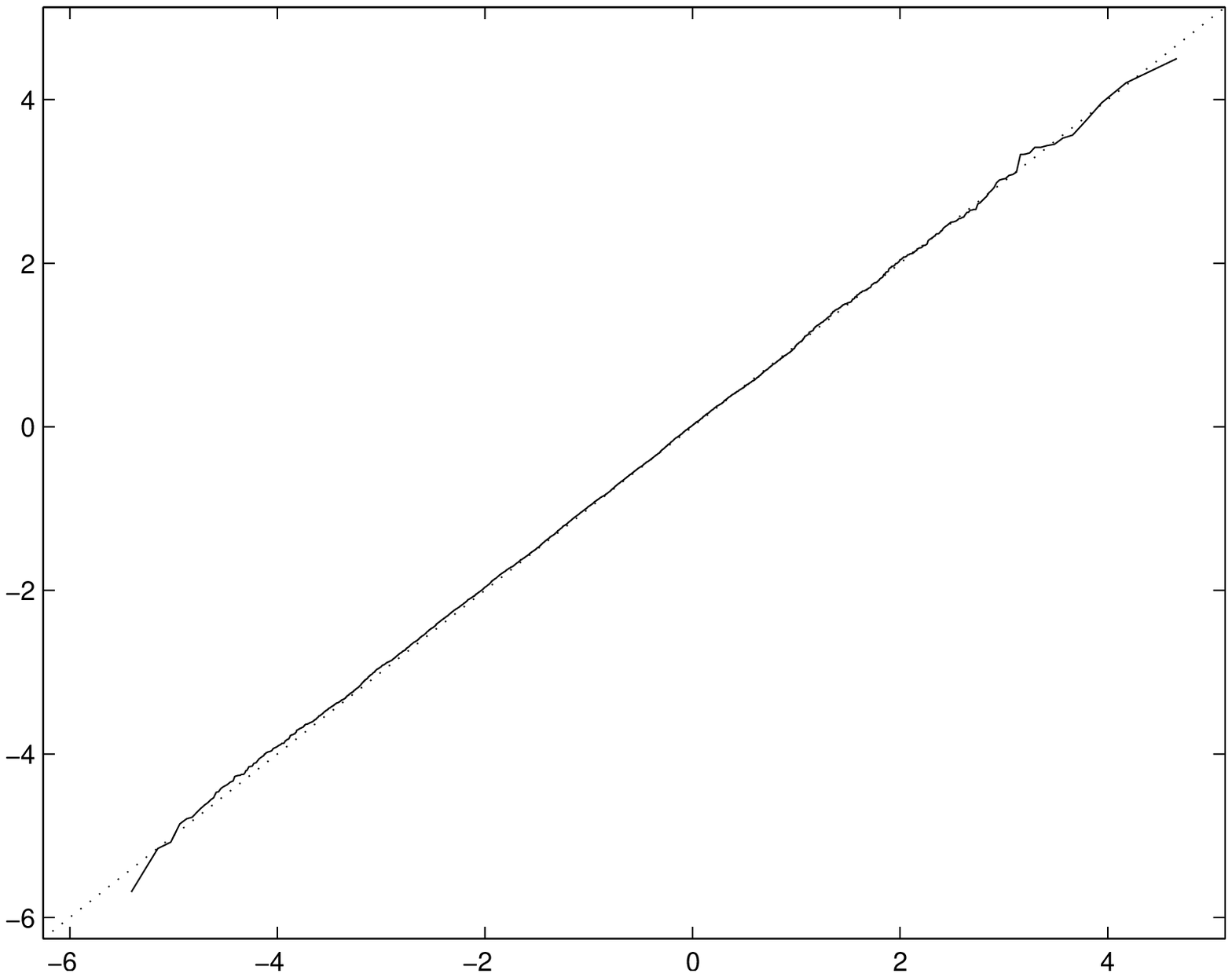} & &
\includegraphics[height=4cm,angle=0]{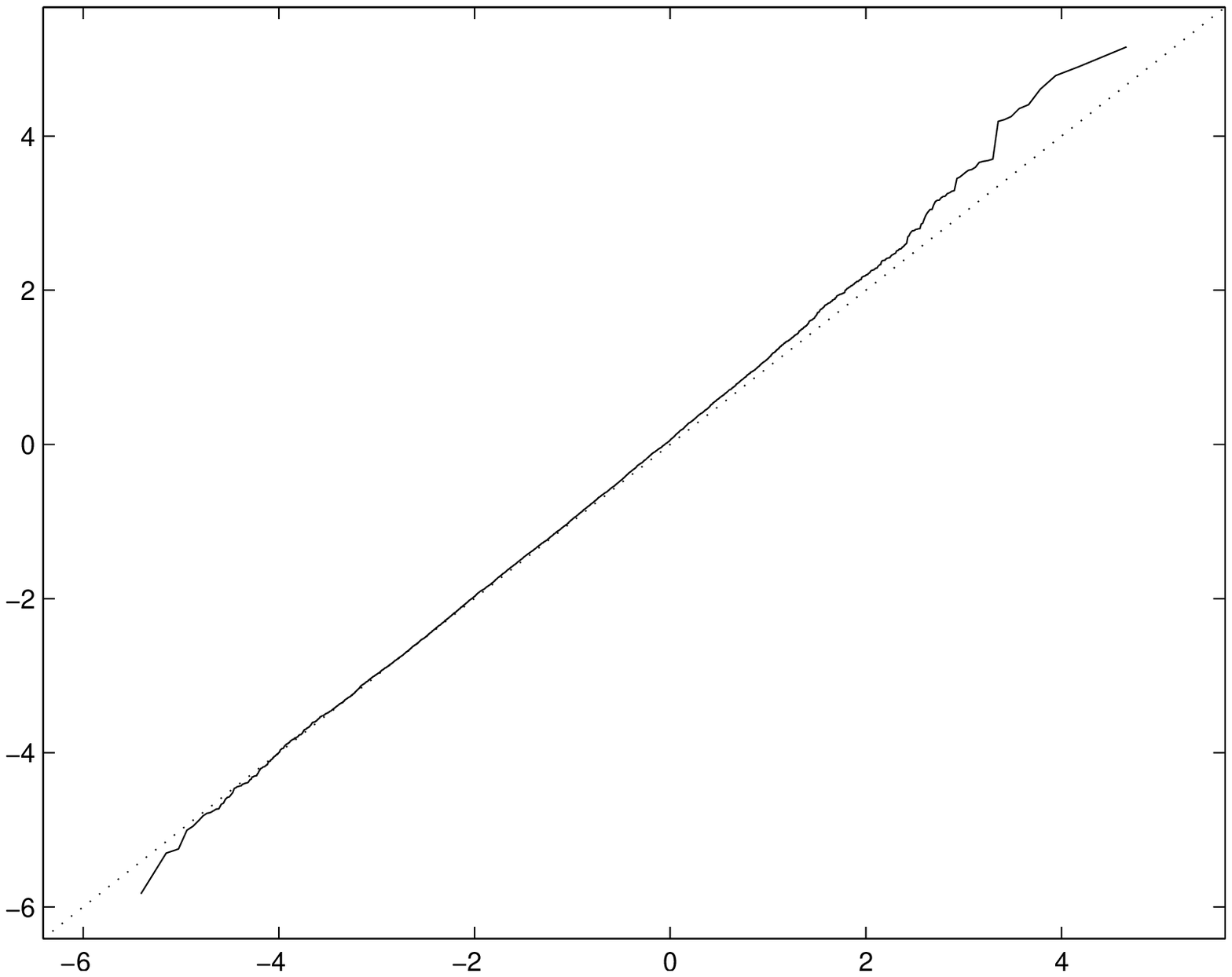} \\
\end{tabular}
\caption{ Probability plot of $10000$ replications of $\lambda_1^{N,p}$ against $F_1^{-1}((i-0.5)/10000)$.
\textit{Left fig.:} Student distribution with 40 degrees of freedom, $N=50$,  $p=200$, $R=10000$. \textit{Right fig.: }Gaussian mixture, $N=50$, $p=200$, $R=10000$.
The largest eigenvalue is rescaled as in (\ref{ressimul})
with $a_1=a_2=0$ (resp. $a_1=a_2=-0.5$) for the Gaussian mixture (resp. Student's) distribution.} 
\label{fig: PPplot}
\end{center}
\end{figure}

\begin{table}\caption{Gaussian mixture, $\gamma_N=50,100$}
\begin{center}\scriptsize
\begin{tabular}{c|ccccccc} 
 $P_c$. & $F_1$ & $10\times 500$ & $25\times 1250$  & $50 \times 2500$ & $10\times 1000$ & $25\times 2500$  & $50 \times 5000$\\
\hline

-3.896	 &.01	&.0147	 &.0114 &.0099	&.0180	 &.0159  &.0113\\
-3.516	 &.025	&.0342	 &.0276 &.0248	&.0377	 &.0322  &.0260\\
-3.180	 &.05	&.0615	 &.0553 &.0488	&.0643	 &.0590  &.0501\\
-2.782	 &.10	&.1156	 &.1051 &.1009	&.1148	 &.1069  &.0977\\
-2.088	 &.25	&.2706	 &.2473 &.2418	&.2663	 &.2552  &.2420\\
-1.269	 &.50	&.5024	 &.4906 &.4839	&.5049	 &.4971  &.4792\\
-0.392   &.75	&.7432	 &.7375 &.7354	&.7450	 &.7411  &.7284\\
0.450	 &.90	&.8899	 &.8887 &.8906	&.8882	 &.8884  &.8893\\
0.979	 &.95	&.9435	 &.9434 &.9448	&.9393	 &.9432  &.9425\\
1.454	 &.975	&.9694	 &.9716 &.9706	&.9698	 &.9698  &.9712\\
2.024	 &.99	&.9867	 &.9876 &.9875	&.9875	 &.9878  &.9883\\
\end{tabular}
\end{center}
{\small The entries are i.i.d. with a Gaussian mixture distribution and $a_1=a_2=0$.}
\end{table}
\begin{table}\caption{Student's $t-$distribution, $\gamma_N=50,100$}
\begin{center}\scriptsize
\begin{tabular}{c|ccccccc} 
 $P_c$. & $F_1$ & $10\times 500$ & $25\times 1250$ & $50 \times 2500$ & $10\times 1000$  & $25\times 2500$  & $50 \times 5000$\\
\hline
-3.896 	&.01	&.0094 &.0093  &.0096 &.0100  & .0099  & .0096	\\
-3.516 	&.025	&.0224	&.0246  &.0221 &.0252  & .0221  & .0229	\\
-3.180 	&.05	&.0470 &.0467  &.0479 &.0503  & .0466  & .0468	\\
-2.782 	&.10	&.1006 &.0898	 &.0961 &.0965  & .0936  & .0962	\\
-2.088 	&.25	&.2415 &.2356  &.2422 &.2466  & .2412  & .2419	\\
-1.269 	&.50	&.4938 &.4819  &.4901 &.4882  & .4893  & .4883	\\
-0.392 	&.75	&.7418 &.7409	 &.7456 &.7418  & .7400  & .7452	\\
0.450 	&.90	&.8970	&.8928  &.8984 &.8953  & .8917  & .8979	\\
0.979 	&.95	&.9468 &.9462  &.9469 &.9495  & .9480  & .9466	\\
1.454 	&.975	&.9717 &.9747	 &.9727 &.9743  & .9733  & .9715	\\
2.024 	&.99	&.9889 &.9889  &.9885 &.9895  & .9902  & .9887	\\
\end{tabular}
\end{center}
{\small The entries are i.i.d. with a $t-$distribution and $a_1=a_2=-0.5$.}
\end{table}

\subsection{Sketch of the proof} 
We here give the main ideas of the proof of both Theorems \ref{theo: unigammafini} and \ref{theo: unigammainfini}.
The proof follows essentially the strategy introduced in
\cite{SosWish} and we refer to this paper for most of the detail. We focus on the case where 
$\gamma=\lim_{N \to \infty}p/N <\infty.$
Basically we compute the leading term in the asymptotic expansion of expectations of traces of
high powers of $M_N$: \be \mathbb{E}\:\left [ \text{Tr}\left (
\frac{1}{N}XX^*\right)^{s_N}\right].\label{largetrace} \ee Here $s_N$ is a sequence such that there exists some constant $c>0$ with $\lim_{N \to \infty}\frac{s_N}{N^{2/3}}=c.$ It is indeed expected that the largest eigenvalues
exhibit fluctuations in the scale $N^{-2/3}$ around $u_+:=\sigma^2(\sqrt{ \gamma_N}+1)^2.$
The core of the proof is to show that for large powers $s_N\sim N^{2/3}$, for any integer $K\geq 1,$ and any real numbers $t_i>0, i=1,\ldots ,K,$ chosen in a compact interval of $\mathbb{R}_+^*$, there exists $C(K)>0$ such that
$\mathbb{E} \Big [\prod_{i=1}^{K}\text{Tr}\left ( \frac{XX^*}{Nu_+}\right)^{[t_iN^{2/3}]}\Big ]\leq C(K)$ and
\begin{equation} 
\Big |\mathbb{E} \Big [\prod_{i=1}^{K}\text{Tr}\left ( \frac{XX^*}{Nu_+}\right)^{[t_iN^{2/3}]}
\Big ]-\mathbb{E}\Big [\prod_{i=1}^{K} \text{Tr}\left ( \frac{X_GX_G^*}{Nu_+}\right)^{[t_iN^{2/3}]}\Big ]\Big |=o(1).\label{unimoments}\end{equation}
Formula (\ref{unimoments}) claims universality of moments of traces of powers of $M_N$ in the scale $N^{2/3}$. Using the machinery developed in \cite{Sos} (Sections 2 and 5) and \cite{SosWish} (Section 2), we can then deduce that the limiting joint distribution of any fixed number of largest
eigenvalues for sample covariance matrices of type $(i)$ to $(i\nu)$ (resp. $(i')$ to ($i\nu'$))
is the same as for complex (resp. real) Wishart ensembles. Here we roughly give the main idea. On one hand, the Laplace transform of the joint distribution of a finite number of the rescaled eigenvalues $\mu_i$ can be conveniently expressed in terms of joint moments of traces as in (\ref{largetrace}). On the other hand, the asymptotic distribution of these rescaled largest eigenvalues (and also the corresponding Laplace transform) is well-known in the Wishart setting. One can then deduce from universality of moments of traces that the asymptotic joint distribution of the largest eigenvalues for any ensemble considered here is the same as for the corresponding Wishart ensemble.
The detail of the derivation of such a result from formula (\ref{unimoments}), including the required asymptotics of correlation functions for Wishart ensembles, can be found in \cite{Sos}, \cite{SosWish} and \cite{NEK}. 
The improvement we obtain with respect to \cite{SosWish} is actually that Formula (\ref{unimoments}) holds for any value $\gamma$. Our result is due to a refinement in the counting procedure of \cite{SosWish}.

\paragraph{}
The paper is organized as follows. In Section \ref{Sec : combiobj}, we introduce the so-called
Narayana numbers. These numbers are the major combinatorial tools needed to adapt the computations of \cite{SosWish} to sample covariance matrices of any sample size to dimension ratio $\gamma_N.$ We also establish a central limit theorem for traces of high powers of $M_N$. Section \ref{Sec: estimatates} is simply a
mimicking of the computations made in \cite{SosWish} and essentially yields formula (\ref{unimoments}). Finally, in Section \ref{Sec: gammainf}, we consider the case where $\gamma_N \to \infty$, which requires some minor modifications.
\paragraph{Acknowledgments}I thank C. Tracy for the Mathematica code of the Tracy-Widom distribution, A. Soshnikov, D. Paul, N. Patterson, R. Cont, L. Choup and C. Semadeni for their great help in the improvement which lead to the final version of this paper. This work was done while visiting UC Davis.

\section{Combinatorics\label{Sec : combiobj}}
In this section, we define the combinatorial objects suitable for the computation of moments of
the spectral measure of random sample covariance matrices.
These combinatorial objects are the Narayana paths and are directly related to the so-called Marchenko-Pastur
distribution. Then, we give the basic technical estimates needed to compute the moments of traces of powers of sample covariance matrices. We illustrate our counting strategy by giving a refinement of the Marchenko-Pastur theorem and also obtain a Central Limit Theorem.
\subsection{Dyck paths and Narayana numbers}
Let $s_N$ be some integer that may depend on $N$. 
Developing (\ref{largetrace}), we obtain
\begin{eqnarray}
&&\mathbb{E} \left [ \text{Tr}\left ( XX^*\right)^{s_N}\right]\crcr
&&
=\sum_{i_o, \ldots ,i_{s_N-1}}\sum_{j_o, \ldots, j_{s_N-1}}
\mathbb{E}\left ( X_{i_o j_o}\overline{ X_{i_1j_o}} \cdots X_{i_{s_N-1}j_{s_N-1}}\overline{X_{i_o j_{s_N-1}}}\right ),
 \label{developtrace}\\
&&\text{where }i_k \in \{1, 2, \ldots, N\} \text{ and }j_k \in \{1, \ldots, p\}, 0\leq k\leq  s_N-1.
\label{ruleR}
\end{eqnarray}
In the whole paper, we denote by $R$ the rule (\ref{ruleR}) for the choice of indices in
(\ref{developtrace}). We shall later prove that such a rule plays a fundamental role in the asymptotics of (\ref{largetrace}).
To each term in the expectation (\ref{developtrace}), we associate three combinatorial objects
that will be needed in the following.

\paragraph{} First, to each term $X_{i_o j_o}\overline{ X_{i_1j_o}} \cdots
X_{i_{s_N-1}j_{s_N-1}}\overline{X_{i_o j_{s_N-1}}}$ occuring in (\ref{developtrace}), we associate the following `` edge path'' $\mathcal P_E$, formed with oriented edges (read from bottom to top)
\be
\begin{pmatrix}j_o\\ i_o \end{pmatrix} \begin{pmatrix}j_o\\ i_1 \end{pmatrix}
\begin{pmatrix}j_1\\ i_1 \end{pmatrix}\cdots  \begin{pmatrix}j_{s_N-1}\\ i_{s_N-1}\end{pmatrix}
\begin{pmatrix}j_{s_N-1}\\  i_o \end{pmatrix}.\label{edgepath}
\ee Due to the symmetry assumption on the entries of $X$, the sole paths leading to a non zero
contribution in (\ref{developtrace}) are such that each oriented edge appears an even number of
times. From now on, we consider only such even edge paths.

\paragraph{} To such an even edge path, we also associate a so-called Dyck path, which is a trajectory $x(t), 0\leq t\leq 2s_N,$ of a simple random walk on the positive half-lattice such that
$$x(0)=0, \: x(2s_N)=0; \: \forall t\in [0, 2s_N], \: x(t)\geq 0\quad \text{and }
x(t)-x(t-1)=\pm 1 .$$
We start the path at the origin and draw up steps $(1, +1)$ and down steps $(1,-1)$ as follows.
We read successively the $2s_N$ edges of (\ref{edgepath}), reading each edge from bottom to top.
Then if the edge (oriented) is read for an odd number of times, we draw an up step. Otherwise we
draw a down step. We obtain in this way a trajectory with $s_N$ up and $s_N$ down steps, which is clearly a Dyck path.
We shall now estimate the number of possible trajectories associated to the edge path. Due to the
constraint on the choices for vertices, we shall distinguish trajectories with respect to the number
of up steps performed at an odd instant. Indeed, they are the moments of time where the vertices can be chosen in the set $\{1, \ldots, p\}$.\\
In the whole paper, we denote by $k$ the number of up steps performed at an odd instant in a Dyck path. 
In particular, we denote by $X_k$ the trajectory associated to $\mathcal{P}_E$, where $k$ is the number of its odd up steps. We also call $\mathcal{X}_{s_N,k}$ the set of Dyck paths of length $2s_N$ with $k$ odd up steps.

\bp \label{Prop: Narayana}\cite{Chen}
Let $\mathbf{N} (s_N, k)$ be the so-called $k$th Narayana number defined by
\be \label{Narayananumber}\mathbf{N}(s_N,k)=\frac{1}{s_N}C_{s_N}^k C_{s_N}^{k-1}.\ee
Then $\mathbf{N}(s_N,k)=\sharp \mathcal{X}_{s_N,k}.$
\ep

\brem For more details about Narayana numbers and their occurrences in various combinatorial problems, we refer the reader to the work of Sulanke \cite{Sulanke1}, \cite{Sulanke2} as well as Stanley \cite{Stanley}. 
\erem

Narayana numbers are intimately linked with Dyck paths. Let $D(2s_N)=\frac{1}{s_N}C_{2s_N}^{s_N+1}$ be the Catalan number counting the number of Dyck paths of length $2s_N.$ It is obvious that
$\sum_{k=1}^{s_N}\mathbf{N} (s_N, k)=D(2s_N).$
Narayana numbers are also linked to the moments of the Marchenko-Pastur distribution defined in (\ref{MarP}), since the following was proved by Jonsson \cite{Jonsson} (see also \cite{Petz} and \cite{BaiMethods}).
\bp \label{Prop: narampdirect1} For any integer $L$, one has that 
\be \lim_{N \to \infty}\frac{1}{N}\mathbb{E}\Bigl [ \emph{Tr} M_N^{L}\Bigr ]=\sigma^{2L}\sum_{k=1}^{L} \gamma^k \mathbf{N}(L, k)= \int x^{L}d\rho_{MP}(x).\label{obsegqmoments}\ee
\ep
 
\brem Proposition \ref{Prop: narampdirect1} was actually proved for a broader class of sample covariance matrices than that considered in this paper. 
\erem

\paragraph{}Last, we associate to the edge path $\mathcal{P}_E$ a ``usual'' path, which we denote by $P_k,$ as
follows.  We mark on the underlying trajectory $x$ the successive vertices met in the edge path.
 The path $ P_k$ associated to (\ref{edgepath}) is then
$i_o\:j_o\:i_1\:j_1\ldots j_{s_N-1}\:i_o.$ For instance, the path associated to
the path $\mathcal{P}_E=$\\ $\begin{pmatrix}1\\ 9 \end{pmatrix}\!
\begin{pmatrix}1\\ 1 \end{pmatrix}\!\begin{pmatrix}3\\ 1 \end{pmatrix}
\!\begin{pmatrix}3\\ 2 \end{pmatrix} \!\begin{pmatrix}4\\ 2 \end{pmatrix}\!\begin{pmatrix}4\\ 1 \end{pmatrix} \!\begin{pmatrix}3\\ 1 \end{pmatrix}\!\begin{pmatrix}3\\ 4 \end{pmatrix}\!\begin{pmatrix}4\\ 4 \end{pmatrix}
\!\begin{pmatrix}4\\ 2 \end{pmatrix}\!\begin{pmatrix}3\\ 2 \end{pmatrix}\!\begin{pmatrix}3\\ 4 \end{pmatrix}\!\begin{pmatrix}4\\ 4 \end{pmatrix}\!\begin{pmatrix}4\\ 1 \end{pmatrix}\!\begin{pmatrix}1\\ 1 \end{pmatrix}\!\begin{pmatrix}1\\ 9 \end{pmatrix}$
is \\given on Fig. \ref{fig: cheminpk} below. 
\vspace*{-0.2cm}
\begin{figure}[htbp]
 \begin{center}
 \begin{tabular}{c}
 \epsfig{figure=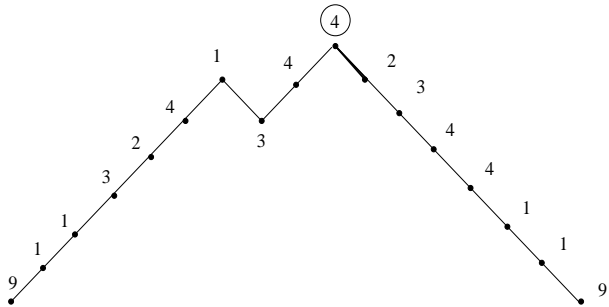,height=4cm,width=8cm, angle=0}
 \end{tabular}
 \caption{ The path $P_k$ with $k=4$, $s_N=8$. \label{fig: cheminpk}}
\vspace*{-0.3cm}
 \end{center}
\end{figure}
\paragraph{}
The three structures $\mathcal{P}_E, X_k$ and $P_k$ we have here introduced will now be used to compute the moments of traces of
(high) powers of $M_N.$ Our counting strategy is as follows. Given the sublying Dyck path $X_k$, we shall estimate the number of edge paths that can be associated to this Dyck path. We shall also estimate their contribution to the expectation (\ref{developtrace}). This is the object
of the two next subsections. 
\subsection{Marked vertices}
In this subsection, we bring out the connection between Narayana paths and the restriction for the choices
 of vertices occuring in the path imposed by the rule $R.$
Given a trajectory $\{x(t), t\leq 2s_N\}\in \mathcal{X}_{s_N,k}$, we shall now count the number of ways to mark
the vertices using the rule $R.$ In this way, we count the number of paths $P_k$ associated to a given trajectory.
The terminology we use is close to the one used in \cite{Sos}, \cite{Si-So1}, \cite{Si-So2} and \cite{SosWish}. We recall the main definitions that will be needed here
and also assume that the reader is acquainted with most of the techniques used in the above papers.

\paragraph{}
The first task is to choose the pairwise distinct vertices occuring in the path. There are at most $p^{s_N+1}$ such vertices.
We shall now define ``marked vertices", separating the cases where they are marked at an odd or even instant.

\bdefi An instant is said to be marked if it is the right endpoint of an up edge of the trajectory $x$.
\edefi

Marked instants correspond to the moments of time where, considering the top and bottom lines separately, one can possibly ``discover'' some vertex not already encountered.  
Consider first the vertices on the top line of the edge path $\mathcal P_E$, that is, vertices occuring at odd instants in the path $P_k$. For $0\leq i \leq s_N$,
call ${\cal T}_i$ the class of vertices of $\{1, \ldots , p\}$ occuring $i$ times as a marked
vertex at an odd instant. Then if we set $p_i=\sharp {\cal T}_i,$ one has 
$$ p=\sum_{i=0}^{s_N}p_i \text{ and } \sum_{i \geq 1}i p_i=k. $$
Note that each time we ``discover" on the top line some new vertex, the corresponding instant is necessarily marked.
Consider also the vertices on the bottom line.
For $0\leq i\leq s_N$, denote by ${\cal N}_i$ the class of vertices of $\{1, \ldots , N\}$ occuring $i$ times as a marked vertex
at an even instant. Then one has, if $n_i= \sharp {\cal N}_i,$
$$ N=\sum_{i=0}^{s_N}n_i \text{ and } \sum_{i \geq 1}i n_i=s_N- k.$$
Note that a vertex from the set $\{1,2, \ldots, N\}$ can occur as a marked vertex on both lines. Yet it
is the type of the vertex on each line which is here taken into account. 
Thanks to the above definition, we characterize a path $P_k$ by its type
$$(n_o, n_1, \ldots , n_{s_N})\:  (p_o, p_1, \ldots , p_{s_N}),
\text{ with } n_i=0, \forall \, i> s_N-k, p_i=0 \:\forall i >k.$$ 
For short, we denote by $(\tilde n, \tilde p)$ the type of such a path.
We also use the following notations. Any vertex $v \in \cup_{i\geq 2}\mathcal{T}_i$ (resp. $v\in \cup_{i\geq 2}\mathcal{N}_i$) is said to be a vertex of self-intersection on the top (resp. bottom) line. A vertex $v\in \mathcal{T}_i$ (resp. $v\in \mathcal{N}_i$) is said to be of type $i$ on the top (resp. bottom) line. 

\paragraph{}
The choice of marked vertices is enough to determine the distinct vertices of the path $P_k$,
if the origin of the path also occurs as a marked vertex. We will see that for typical paths, this is
not the case and $i_o \in \mathcal{N}_0.$
Thus, given the type $(\tilde n, \tilde p)$ of the path, the number of ways to assign vertices at the marked instants and choose the origin is then at most:
\be \label{choixdessommets}N \frac{N!}{\prod_{i=0}^{s_N}n_i!} \frac{p!}{\prod_{i=0}^{s_N} p_i!}
 \frac{k!}{\prod_{i\geq 2}(i!)^{p_i}}\frac{(s_N-k)!}{\prod_{i\geq 2}(i!)^{n_i}}.
\ee
Indeed, one distributes the vertices of $\{1, \ldots, N\}$ and $\{1, \ldots, p\}$ into the possible classes $\mathcal{N}_i, \mathcal{T}_i, 1\leq i\leq s_N,$ choose the corresponding marked occurrences of each vertex and fix the origin.
Once the marked vertices and the origin of the path are chosen, there remains to fill in
the blanks of the path $P_k$. Due to self-intersections, there are multiple ways to do so.
We investigate this numbering in the sequel and consider at the same time the expectation of the filled path. 
\subsection{Filling in the blanks of the path}
Consider now a path with $k$ odd marked instants and of type $(\tilde n, \tilde p)$.
Call $\Omega _{max}$ the maximum number of ways to fill in the blanks of $P_k$ at the
unmarked instants, once the marked instants and the origin are given.

\bp \label{Prop :closingandweight}Set $\Omega _{max} \mathbb{E}_{max}:=\max_{(\tilde n, \tilde p)}
\Omega _{max}\Big |\mathbb{E} \Big [\prod_{j =0}^{s_N-1} M_{i_ji_{j+1}}\Big ]\Big | $.\\ There exists $\tilde C>0$
independent of $p,N, k$ and $s_N$ such that
\begin{equation}\label{majoWmE}
 \Omega _{max} \mathbb{E}_{max}  \leq 2
\frac{\sigma^{2s_N}}{N^{s_N}}\prod_{l= 2}^{s_N-k}\left (\tilde Cl\right)^{ln_l}\prod_{m= 2}^{k}\left (\tilde Cm\right)^{mp_m}.
\end{equation}
\ep

\paragraph{Proof of Proposition \ref{Prop :closingandweight}}We only sketch the proof
which follows essentially the same steps as that of Lemma 1 in \cite{Si-So2}.
Assume that in $P_k$, at the unmarked instant $t$, one makes a down step with left vertex $i$. If $i$ is of type $1$, then there is no choice for the right endpoint of such an edge. In general, the maximal number of possible right endpoints depends on the multiplicity of $i$ as a marked vertex.
Knowing the parity of $t$, one also knows whether $i$ is a vertex on the top or on the bottom
line of $\mathcal{P}_E$.
Thus the sole top or bottom multiplicity of the vertex has to be taken into account to estimate
the number of ways to close the edges. Thus, it is not hard to see that the number of ways to close the path is at most
$$ 2\prod_{l= 2}^{s_N-k}\left (2l\right)^{ln_l}\prod_{m= 2}^{k}\left (2m\right)^{mp_m}.$$
The extra factor $2$ comes from the case (negligible) where the origin is of type $1$. To consider simultaneously, the expectation of the path, one also has to take into account the number of times each oriented edge is read.
Assume that an edge $\begin{pmatrix}v\\ w \end{pmatrix}$ is read $2q$ times, with $q\geq 2$. Call $l_u(v;w)$ (resp. $l_d(w;v)$) the number of times $v$ (resp. $w$) is a marked vertex of this edge. Then, if $l_d(w;v)l_u(v;w)>0$, one has that
$\mathbb{E}|M_{vw}|^{2q}\leq (\tau q)^q \leq \left (2\tau l_u(v;w)\right)^{l_u(v;w)} \left (2\tau l_d(w;v)\right)^{l_d(w;v)}.$
Now if an oriented edge is read $2l$ times then it is closed $l$ times along the same edge. That is, we overcount the number of ways to close the path.
Thus, if we let $2l(ij)$ be the number of times the oriented edge $(ij)$ is read in the path, we obtain that
\begin{eqnarray*}
&\dfrac{\Omega _{max} \mathbb{E}_{max}}{2\sigma^{2s_N}}&\leq  
\prod_{l(ij)>1} \frac{\left (C_{\tau} l(ij)\right)^{l(ij)}}{l(ij)!}
\prod_{l= 2}^{s_N-k}\left (2l\right)^{ln_l}\prod_{m= 2}^{k}\left (2m\right)^{mp_m}\crcr
&&\leq \prod_{l= 2}^{s_N-k}\left (\tilde Cl\right)^{ln_l}\prod_{m= 2}^{k}\left ( \tilde Cm\right)^{mp_m},
\end{eqnarray*}
where $C_{\tau},\tilde C$ are some constants independent of $p,$ $k,$ $N$ and $s_N.$ $\square$
\brem \label{Rem: cutoff} In the case where the entries $X_{ij}$ have polynomial tails, with $P(|X_{ij}|\geq x)\leq (1+x)^{m_o}$ for some $m_o>36$, one can first consider (up to a set of negligible probability) that all the entries of $X/\sigma$ are smaller in absolute value than $\Gamma_N:=N^{2/m_o+\epsilon}$ for some $\epsilon>0$ small enough. This is true if $\gamma<\infty.$ Then Proposition \ref{Prop :closingandweight} has to be replaced with  
$$ \frac{\Omega _{max} \mathbb{E}_{max}}{2\sigma^{2s_N}N^{-s_N}}  \leq \prod_{l= 2}^{s_N-k}\left (lC(1+\Gamma_N^4 1_{l\geq \frac{m_o}{4}})\right)^{ln_l}\prod_{m= 2}^{k}\left (mC(1+\Gamma_N^4 1_{m\geq \frac{m_o}{4}})\right)^{mp_m}.$$
This follows from the fact that to each edge seen $2l\geq 2m_o$ times there corresponds at least $l/2$ marked occurences of one of its endpoints. 
\erem
\paragraph{}In the two following subsections, we investigate moments of Traces of powers of $M_N$ in scales $s_N<<\sqrt N.$ This will give the foundations for the asymptotics of higher moments.
\subsection{Narayana numbers and the Marchenko-Pastur distribution}
In this subsection, we illustrate our counting strategy and present a refinement of (\ref{obsegqmoments}), which allows to consider higher moments than in Proposition \ref{Prop: narampdirect1}. 

\bp \label{Prop: narampdirect}If $s_N <<\sqrt{N},$ one has that 
$$\frac{1}{N}\mathbb{E}\Bigl [ \emph{Tr} M_N^{s_N}\Bigr ]=\sigma^{2s_N} \sum_{k=1}^{s_N} \gamma_N^{k}
\mathbf{N}(s_N,k)(1+o(1)).$$
\ep

\paragraph{Proof of Proposition \ref{Prop: narampdirect}}
The proof is similar to that of the classical Wigner theorem using Dyck paths (see e.g. \cite{BaiMethods}). It is divided into two steps. First, we show
that paths for which $ \sum_{i\geq 2}n_i+p_i>0$ yield a negligible contribution to $\mathbb{E}\left [\text{Tr} M_N^{s_N}\right]$.
Then we estimate the contribution of paths with vertices of type 1 at most, which give the leading term of the asymptotic expansion of $\mathbb{E}\left [\text{Tr} M_N^{s_N}\right]$, as long as $s_N<<\sqrt N.$

\paragraph{}
Denote by $Z(k,(\tilde n, \tilde p))$ the contribution of paths with $k$ odd marked instants and of type
$(\tilde n, \tilde p).$
Using Proposition \ref{Prop :closingandweight} and (\ref{choixdessommets}), we deduce that
\begin{eqnarray}
&&Z(k,(\tilde n, \tilde p))\crcr
&&\leq \mathbf{N} (s_N, k)\frac{2\sigma^{2s_N}}{N^{s_N}}N\frac{p!}{p_o!}\frac{N! }{n_o!} \frac{k!}{p_1!}\frac{(s_N-k)!}{n_1!}
\prod_{i=2}^{s_N}\frac{\left (\tilde Ci\right )^{in_i+ip_i}}{p_i! n_i!(i!)^{p_i} (i!)^{n_i}} \crcr
&&\leq \mathbf{N} (s_N, k) N2\sigma^{2s_N}\prod_{i=2}^{s_N}\frac{\left (Ck\right )^{ip_i} \,\left( C(s_N-k)\right) ^{ in_i}}{p_i!n_i!}\:\dfrac{ N^{\sum_{i\geq 1}n_i}\: p^{\sum_{i\geq 1}p_i}}{N^{\sum_{i\geq 1} in_i+ip_i}}\:  
\crcr
&&\leq \mathbf{N} (s_N, k) N 2\sigma^{2s_N}\gamma_N^k \prod_{i\geq 2}\frac{1}{n_i!}\left( \frac{C^i(s_N-k)^i}{N^{i-1}}\right)^{n_i}
\prod_{i\geq 2}\frac{1}{p_i!}\left( \frac{C^ik^i}{p^{i-1}}\right)^{p_i}, \label{contrigal}
\end{eqnarray}
where in the last line we have used that $\gamma_N^{\sum_{i\geq 1}p_i}=\gamma_N^{k-\sum_{i\geq 2}(i-1)p_i}$ and $C>0$ is a constant independent of $N,p,k$ and $s_N$ (whose value may change from line to line).

\paragraph{}We denote by $Z_2$ the contribution of paths for which $ \sum_{i\geq 2}n_i+p_i>0$.
By Proposition \ref{Prop :closingandweight}, and using summation, one has that
$$
\frac{Z_2}{2\sigma^{2s_N}}\leq \sum_{k=1}^{s_N} N \mathbf{N}(s_N,k)\gamma_N^k \sum_{\tilde M_1, \tilde M_2:\tilde M_1+ \tilde M_2>0}
\frac{\left ( \frac{2Cs_N^2}{N}\right)^{\tilde M_1}
\left ( \frac{2Cs_N^2}{p}\right)^{\tilde M_2}}{\tilde M_1!\tilde M_2!},
$$
where $\tilde M_1=\sum_{i\geq 2}n_i$ and $\tilde M_2=\sum_{i\geq 2}p_i$.
Thus, it is straightforward to see that there exists some constant $B>0$ independent of $N$ such that
\be \label{estz2}\frac{Z_2}{N}\leq B\frac{s_N^2}{N}\times \sigma^{2s_N} \sum_{k=1}^{s_N} \gamma_N^{k}
\mathbf{N}(s_N,k).\ee 
From this, we can deduce that $Z_2/N=o(\sigma^{2s_N}(1+\sqrt{\gamma_N})^{2s_N})$. 

\paragraph{}
We now show that only paths with vertices of type 1 (except the origin which is unmarked) have to
be taken into account. In this case, once the vertices occuring in the path have been chosen, there is no choice for filling in the blanks of the path $P_k.$
Furthermore, each edge is passed only twice in the path $\mathcal{P}_E$, once at an odd instant and once at an even instant. Thus, denoting by $Z_1:=\sum_{k=1}^{s_N}Z_1(k)$ the contribution of such paths, one has that 
$$Z_1= \sum_{k=1}^{s_N}N \mathbf{N}(s_N,k)\gamma_N^k \sigma^{2s_N}\prod_{i=1}^{s_N}\frac{N-i}{N}=(1+o(1)) \sum_{k=1}^{s_N}N \mathbf{N}(s_N,k)\gamma_N^k \sigma^{2s_N}.$$
Using (\ref{estz2}), this finishes the proof that $Z_2=o(1)Z_1$. The contribution of paths with marked origin and vertices of type $1$ at most is of order $Z_1s_N/N$ and thus negligible.
This finishes the proof of Proposition \ref{Prop: narampdirect}. $\square$

\brem \label{rem : u+} Set $\hat k=[\frac{\sqrt{\gamma_N}}{1+\sqrt{\gamma_N}}s_N]+1.$ 
For any sequence $s_N>>1$, one has that $ \sum_{1\leq k \leq s_N} \frac{N}{s_N}C_{s_N}^{k}C_{s_N}^{k-1} \gamma_N^{k}\sigma^{2s_N}=O(u_+^{s_N})$ and that the main contribution to the expectation (\ref{largetrace}) should come from paths with $\hat k(1+o(1))$ odd marked instants.
Indeed, using Stirling's formula, one has
$$\max_{1\leq k \leq s_N}\left (C_{s_N}^k\right)^2 \gamma_N^k \sim
\left (C_{s_N}^{\hat k}\right)^2\gamma_N^{\hat k}\sim \left ( 1+\sqrt{\gamma_N}\right)^{2s_N}\frac{(1+\sqrt{ \gamma_N})^2}{s_N \sqrt{\gamma_N}}\frac{1}{2\pi}.$$
It is also easy to check that, for any $l>0$, one has that
\begin{equation}
\mathbf{N}(s_N,\hat k+l)\gamma_N^{\hat k+l}\leq \mathbf{N}(s_N,\hat k)\gamma_N^{\hat k} \exp{\{-\frac{C_{\gamma} l^2}{(s_N-\hat k)}\}}, \label{decroissancel2}
\end{equation}
for some constant $C_{\gamma}$ depending on $\gamma$ only. In the case where $l<0$, we fix some $\Delta >0$ large. Then one can show that, for any $-\Delta(s_N-\hat k+1)<l<0$, 
$\frac{\mathbf{N}(s_N,\hat k-1+l)\gamma_N^{\hat k-1+l}}{\mathbf{N}(s_N,\hat k-1)\gamma_N^{\hat k-1}}
\leq \exp{\{- \frac{l^2}{(2\Delta +2)(s_N-\hat k+1)}\}}.$ One also has that for any $l\leq -\Delta(s_N-\hat k+1)$, 
$\frac{\mathbf{N}(s_N,\hat k-1+l)\gamma_N^{\hat k-1+l}}{\mathbf{N}(s_N,\hat k-1)\gamma_N^{\hat k-1}}\leq e^{-\Delta(s_N-\hat k)/3}.$
This ensures that $ \sum_{1\leq k \leq s_N} \frac{N}{s_N}C_{s_N}^{k}C_{s_N}^{k-1} \gamma_N^{k}\sigma^{2s_N}=O(u_+^{s_N}),$
yielding Remark \ref{rem : u+}. \erem
\brem \label{Rem: cutoffmaj} In the case of polynomial tails, (\ref{contrigal}) is multiplied by a factor 
$\prod_{i \geq m_o/4} \Gamma_N^{4i(n_i+p_i)}.$ If $\epsilon <\frac{m_o-36}{12m_o}$ and $s_N=O(N^{2/3}),$ (which is the largest scale considered in this paper), this has no impact on computations as $\frac{(s_N\Gamma_N^{4})^i}{N^{i-1}}<<1$ for any $i\geq m_o/4$. All the results stated in the following can be proved in the case of polynomial tails up to minor technical modifications (which amounts essentially to considering apart vertices of type at least $m_o/4$).
\erem
\subsection{\label{subsec : clt}A Central Limit Theorem}
The main result of this subsection is the following Proposition. Set $u_+=\sigma^{2}(1+\sqrt{\gamma_N})^2.$ We show that all the moments of $\text{Tr}(M_N/u_+)^{s_N}$ are bounded and universal, as long as $1<<s_N<<\sqrt{N}$. Assume that $\lim_{N \to \infty}\gamma_N=\gamma<\infty$ and set  $l_{\beta}=1/(\beta \pi)$ where $\beta=1$ (resp. $\beta=2$) in the case where $M_N$ is real (resp. complex). 
\bp \label{Prop: variance} Assume that $1<<s_N<<\sqrt{N}$ and set $\tilde M_N=\frac{M_N}{u_+}.$ Then, there exists $D>0$ such that $\text{Var}\left (\emph{ Tr} {\tilde M_N}^{s_N}\right)\leq D,$ for any $N$,  
and $\displaystyle{\lim_{N \to \infty}\text{Var}\left (\emph{ Tr} {\tilde M_N}^{s_N}\right)=l_{\beta}}.$
Similarly, for any integer $k$,  
\begin{eqnarray*}
&&\mathbb{E}\left [ \emph{Tr}{\tilde M_N}^{s_N}- \mathbb{E}[\emph{Tr}{\tilde M_N}^{s_N}]\right]^{2k}=(2k-1)!!\: l_{\beta}^k(1+o(1)),\crcr
&&\mathbb{E}\left [ \emph{Tr}{\tilde M_N}^{s_N}- \mathbb{E}[\emph{Tr} {\tilde M_N}^{s_N}]\right]^{2k+1}=o(1).
\end{eqnarray*}
\ep

\brem \label{rem: normal} In \cite{Jonsson}, a Central Limit Theorem (CLT) is also established for traces of fixed (independent of $N$) moments of $M_N$. In this case, the limiting Gaussian distribution does depend on the fourth moment of the law of the entries. The above CLT is also stated but not proved in Remark 6 of \cite{Si-So1} (a factor $1/\beta$ is missing) in the case where $\gamma=1$.\erem 

\paragraph{Proof of Proposition \ref{Prop: variance}}We only give the proof for the variance. 
The proof of the asymptotics for higher moments is a rewriting of pages 128-129 in \cite{Si-So2} (also \cite{PecheFeral} Section 6) and is skipped.
In the following, $C_1, \ldots, C_6, C'_2,C'_3$ denote some positive constants independent of $N$.
One has that 
\begin{eqnarray*}
&\text{Var (Tr }M_N^{s_N})=&\frac{1}{N^{2s_N}}\sum_{\mathcal{P}_E,\mathcal{P}'_E}\Big [\mathbb{E}\left (\prod_{e_i\in \mathcal{P}_E } \prod_{e'_i\in \mathcal{P}'_E }\hat X_{e_i}\hat X_{e'_i}\right)\crcr
&&-\mathbb{E}\left (\prod_{e_i\in \mathcal{P}_E} \hat X_{e_i}\right) \mathbb{E}\left ( \prod_{e'_i\in \mathcal{P}'_E}\hat X_{e'_i}\right)\Big ].
\end{eqnarray*}
Here, given an edge $e=(v_1, v_2)$, $\hat X_e$ stands for $X_{v_1v_2}$ if $e$ occurs at an odd instant of $\mathcal{P}_E$ or for $\overline{X_{v_2v_1}}$ if it occurs at an even instant.
Now it is clear that the non zero terms in the above sum come from pairs of paths $\mathcal P_E$, $\mathcal P'_E$
sharing at least one oriented edge and such that each edge appears an even number of times in the union of the two paths. We say that such paths are correlated.
To estimate the number of correlated paths and their contribution to the variance, we use the \emph{construction procedure }defined in Section 3 of \cite{Si-So1}. This construction associates a path of length $4s_N-2$ to a pair of correlated paths.\\
Let $\mathcal P_1$ and $\mathcal P_2$ be two correlated paths of length $2s_N$. When reading the edges of $\mathcal{P}_1$, let $e$ denote the first oriented edge common to the two paths. Let also $t_e$ and $t'_e$ be the instants of the first occurrence of this edge in $\mathcal{P}_1$ and $\mathcal{P}_2$. Then we are going to glue the two paths $\mathcal{P}_1$ and $\mathcal{P}_2$, in such a way that we erase the two first occurrences of $e$ in each of these paths. The glued path, denoted $\mathcal P_1\vee \mathcal P_2$, is obtained as follows.
We first read $\mathcal P_1$ until we meet the left endpoint of $e$ at the instant $t_e$. Then we switch to $\mathcal P_2$ as follows. Assume first that $t_e$ and $t'_e$ are of the same parity. We then read the path $\mathcal P_2$, starting from $t'_e$, in the reverse direction to the origin and restart from the end of $\mathcal{P}_2$ until we come back to the instant $t'_e+1$. If $t_e$ and $t'_e$ are not of the same parity, we read the edges of $\mathcal{P}_2$ in the usual direction starting from $t_e'+1$ and until we come back to the instant $t'_e$. We have then read all the edges of $\mathcal P_2$ except the edge $e$ occuring between $t'_e$ and $t'_{e+1}.$ We then read the end of $\mathcal P_1$, starting from $t_e+1.$
Having done so, we obtain a path $\mathcal P_1 \vee \mathcal P_2$ which is of length $4s_N-2$.
One can also note that the trajectory of $\mathcal P_1 \vee \mathcal P_2$ does not descent lower than the level $x(t_e)$ during the time interval $[t_e, t_e+2s_N-1]$, by the definition of $e$ and $t_e$.\\
Now, to reconstruct the paths $\mathcal P_1$ and $\mathcal P_2$ from $\mathcal P_1 \vee \mathcal P_2$, it is enough to determine the instant at which one has switched from one path to the other, the origin of the path $\mathcal P_2$ and the direction in which $\mathcal P_2$ is read. There are at most $4s_N$ ways to determine the origin and the direction once the instant of switch is known.
To estimate the number of preimages of a given path $\mathcal P_1 \vee \mathcal P_2$ of length $4s_N-2$ and with $k$ odd up steps, one has to give an upper bound for the number of instants $t_e$ in $\mathcal P_1 \vee \mathcal P_2$, which can be the instants of switch. 
To this aim, fix some $t_e \in [0, 2s_N-1]$ and assume that the trajectory of $\mathcal P_1 \vee \mathcal P_2$ does not go below the level $x(t_e)$ during an interval of time of length greater than or equal to $2s_N-1.$ Assume that $x(t_e)>0.$ Set then 
$$l=\inf \{t \geq t_e, x(t)=x(t_e), x(t+1)=x(t_e)-1\}-2s_N+1.$$
Denote by $T_2$ the sub-trajectory in the interval $[t_e, t_e+2s_N-1+l]$. It is a Dyck path. Denote also by $T_1$ the remaining part of the trajectory: it is also a Dyck path, along which the instant $t_e$ has been chosen. We denote by $k_1$ the number of the odd up steps of $T_1$. As the trajectory of $\mathcal{P}_1\vee \mathcal{P}_2$ is obtained by inserting $T_2$ at the instant $t_e$ in $T_1$, and using the fact that $\mathcal{P}_1\vee \mathcal{P}_2$ and $\mathcal{P}_1 \cup \mathcal{P}_2$ have all the same edges but one, one can then deduce (see \cite{Si-So1}, p. 11-13, for the detail) that the contribution of correlated pairs is at most of order  
\begin{eqnarray}
&&\!\!\!\!\!\!\!\!\!\!
\sum_{k=1}^{2s_N-1}\sum_{l=0}^{2s_N-1}\sum_{k_1\leq k \wedge 2s_N-1-l}\frac{\mathbf{N}( s_N-\frac{(1+l)}{2},k_1)\mathbf{N}( s_N+\frac{l-1}{2},k-k_1)}{\mathbf{N}( 2s_N-1,k)}\crcr
&&\!\!\!\!\!\!\!\!\!\!\times (2s_N-1-l)\frac{4s_N\sigma^2}{N} Z_1(4s_N-2,k) +
\crcr
&&\label{contrcorrel}\\
&&\!\!\!\!\!\!\!\!\!\!
\sum_{k=1}^{2s_N-1}\sum_{l=0}^{2s_N-1}\sum_{k_1\leq k \wedge 2s_N-1-l}\frac{\mathbf{N}( s_N-\frac{(1+l)}{2},k_1)\mathbf{N}( s_N+\frac{l-1}{2},s_N+\frac{l-1}{2}+k_1-k)}{\mathbf{N}( 2s_N-1,k)}\crcr
&&\!\!\!\!\!\!\!\!\!\!
\times (2s_N-1-l)\frac{4s_N\sigma^2}{N} Z_1(4s_N-2,k).\crcr &&\label{contrcorrelodd}
\end{eqnarray}
Here $Z_1(4s_N-2,k)$ is the contribution of paths of length $4s_N-2$ with $k$ odd up steps to the expectation 
$\mathbb{E}[ \text{Tr}M_N^{2s_N-1}]$ and $(\ref{contrcorrel})$ (resp. $(\ref{contrcorrelodd})$) corresponds to the case where $t_e$ is even (resp. odd). The term $(2s_N-l-1)$ in (\ref{contrcorrel}) comes from the determination of $t_e$ and where $\sigma^2/N=\mathbb{E}(|X_e|^2)/N$, if $e$ is the edge erased from $\mathcal{P}_1\vee \mathcal{P}_2$. It can indeed be shown that paths for which such an edge occurs also in $\mathcal{P}_1\vee \mathcal{P}_2$ yield a contribution of order $s_N/N$ that of typical paths and are thus negligible.\\
We first show that the variance is bounded. In the following, we set $s_1(l)=s_N-\frac{1+l}{2}$ and $s_2(l)=s_N+\frac{l-1}{2}.$
Considering for instance $(\ref{contrcorrel})$, ($(\ref{contrcorrelodd})$ is similar), it is enough to prove that there exists a constant $C_1>0$ such that 
$$\sum_{l=0}^{2s_N-1}\sum_{k_1\leq k \wedge 2s_N-1-l}\frac{\mathbf{N}( s_1(l),k_1)\mathbf{N}( s_2(l), k-k_1)}{\mathbf{N}( 2s_N-1,k)}(2s_N-1-l)\leq C_1\sqrt{s_N}.$$
One can easily see that it is enough to consider the case where $2s_N-1-l\geq \sqrt{s_N}.$  
It is also straightforward by Remark \ref{rem : u+} and Proposition \ref{Prop: narampdirect} to see that one can choose $0<2\beta'<\frac{\sqrt \gamma}{1+\sqrt \gamma}<2\beta''<2$ such that 
$$\sum_{k\leq \beta' s_N\text{ or }k\geq \beta''s_N }Z_1(4s_N-2,k)<<s_N^{-5}N u_+^{2s_N-1}.$$ 
This is enough to ensure that the contribution of correlated pairs such that the corresponding glued path has $k$ odd up instants for some $k \leq \beta's_N$ or $k\geq \beta'' s_N$ is negligible in the large-$N$-limit.
We now set
\be \label{majoinfini}f(k_1):=\frac{\mathbf{N}( s_1(l),k_1)\mathbf{N}( s_2(l),k-k_1)}{\mathbf{N}( 2s_N-1,k)}.\ee
Then, $l$ and $k$ being fixed, $f$ is maximal at $\tilde k_1= [k \frac{2s_N-1-l}{4s_N-2}](+1).$ Furthermore, one can check that there exist constants $C_2, C_2'>0$ such that $f(\tilde k_1+j)\leq C_2\exp{\{-C_2'j^2/\tilde k_1\}}$ for any $j$. 
From this we deduce that 
\begin{eqnarray}
&&\!\!\!\!\!\!\!\sum_{k_1\leq k \wedge 2s_N-1-l}\frac{\mathbf{N}( s_1(l),k_1)\mathbf{N}( s_2(l),k-k_1)}{\mathbf{N}( 2s_N-1,k)}(2s_N-1-l)\crcr
&&\!\!\!\!\!\!\!
\leq C'_3(2s_N-1-l)^{3/2}\frac{\mathbf{N}( s_1(l),\tilde k_1)\mathbf{N}( s_2(l),k-\tilde k_1)}{\mathbf{N}( 2s_N-1,k)}.\label{majoenko}
\end{eqnarray}
It is now an easy consequence of Stirling's formula that
\be
\sum_{l=0}^{2s_N-1}C_3(2s_N-1-l)^{3/2}\frac{\mathbf{N}( s_1(l),\tilde k_1)\mathbf{N}( s_2(l),k-\tilde k_1)}{\mathbf{N}( 2s_N-1,k)}
\leq C_4 \frac{\sqrt{s_N}}{(1-\alpha_N)^2},\label{majosoml}\ee
where $\alpha_N=k/(2s_N-1).$
Using Proposition \ref{Prop: narampdirect}, one can also show that there exists $C_5>0$ such that $$\sum_{\beta' s_N\leq k\leq \beta'' s_N}\frac{1}{\alpha_N(1-\alpha_N)}Z_1(4s_N-2,k)\sim C_5s_N^{-3/2}N u_+^{2s_N-1}.$$
Combining the whole yields that there exists a constant $C_6$ such that $(\ref{contrcorrel})+(\ref{contrcorrelodd}) \leq C_6u_+^{2s_N}.$ \\
In the case where $x(t_e)=0$, $t_e$ is chosen amongst the returns to the level $0$ of the trajectory. It can be shown that the number of such instants is negligible with respect to $\sqrt{s_N}$ in typical paths. This follows from arguments already used above and in \cite{Si-So1} p. 13.\\ 
To compute the variance, we notice that in (\ref{contrcorrel}), the term $(2s_N-1-l)$ can actually be replaced with $s_1(l)-k_1$. Indeed as $t_e$ is even, the first step after $t_e+2s_N-1+l$ is a down step occuring at an odd instant. Also, there are only $s_N$ choices for the origin of ${\mathcal P}_2$, since one knows the parity of $e$  in $\mathcal{P}_2$ once the orientation of $\mathcal{P}_2$ is fixed.
Then, using (\ref{contrcorrel}), (\ref{contrcorrelodd}), Remark \ref{rem : u+}, the exponential decay of $f(k_1)$, and Proposition \ref{Prop: narampdirect}, one can deduce that (for the real case)
\begin{eqnarray}
&&\lim_{N \to \infty}\text{Var Tr}\tilde M_N^{s_N}\crcr
&&=\lim_{N \to \infty}2s_N(1+\frac{1}{\sqrt \gamma})\frac{1}{1+\sqrt \gamma}\sum_{l\leq 2s_N-1}\frac{s_1(l)}{(1+\sqrt{\gamma})^{4s_N}}\crcr
&&\times \sum_{k\geq 1}\sum_{ k_1\leq k}\gamma_N^{k}\mathbf{N}( s_1(l), k_1)\mathbf{N}( s_2(l),k- k_1)\label{serie}\\
&&=\lim_{N \to \infty}2s_N(1+\frac{1}{\sqrt \gamma})\frac{1}{1+\sqrt \gamma}\sum_{l\leq 2s_N-1}\frac{s_1(l)}{(1+\sqrt{\gamma})^{4s_N}} E_{s_1(l)}E_{s_2(l)}\label{serie2}\\
&&=l_{1}.\nonumber \end{eqnarray}
In (\ref{serie2}),  we have set $E_{k}=\int x^{k}\frac{\sqrt{((1+\sqrt \gamma)^2-x)(x-(1-\sqrt \gamma)^2)}}{2\pi x}dx$, and the equality follows from the fact that $(\ref{serie})$ is a Cauchy product.
The value of $l_{1}$ can be deduced from Formulas 4.7 in \cite{SosWish} and 3.6 in \cite{Si-So1}. The computation of $l_2$ follows from the fact that, in the complex case, the occurences $e$ in $\mathcal{P}_1$ and $\mathcal{P}_2$ cannot have the same parity (in typical paths).
$\square$

\section{The case where $\mathbf{\gamma_N \to \gamma,}$ $\mathbf{1\leq \gamma<\infty}$\label{Sec: estimatates}}
The aim of this section is to prove the following universality results.
Let $M_N=\frac{1}{N}XX^*$ be a sequence of sample covariance matrices satisfying $(i)$ to $(i\nu)$
(resp. $(i')$ to $(i\nu'))$. 
Let $K$ be some given integer and $c_1, \ldots , c_K$ be constants chosen in some compact interval $J\subset \mathbb{R}^+$ (independent of $N$). 
Consider sequences $s_N^{(i)}, i=1, \ldots, K,$ such that $\lim_{N \to \infty}\frac{s_N^{(i)}}{N^{2/3}}=c_i.$ 

\bt \label{theo: unimoments} Set $u_+=\sigma^2(1+\sqrt{\gamma_N})^2.$ There exists a constant $\tilde C_1=\tilde C_1(K, J)>0$ such
that
\begin{eqnarray*}&&\mathbb{E} \left  [ \prod_{i=1}^K \emph{Tr}  \left(\frac{XX^*}{Nu_+} \right)^{s_N^{(i)}}\right]\leq \tilde C_1\text{ and }\crcr 
&&\mathbb{E}\left[ \prod_{i=1}^K \emph{Tr}  \left(\frac{XX^*}{Nu_+}\right)^{s_N^{(i)}}
\right]=\mathbb{E}\left[ \prod_{i=1}^K \emph{Tr}  \left(\frac{X_GX_G^*}{Nu_+}\right)^{s_N^{(i)}}\right](1+o(1)).
\end{eqnarray*} \et

The proof of Theorem \ref{theo: unimoments} is the object of this section. We actually focus on the case where $K=1$. Indeed, the proof of Theorem \ref{theo: unimoments} for $K>1$ is a rewriting of the arguments used in \cite{Sos} (p. 41), Subsection \ref{subsec : clt} and of the arguments used in the case where $K=1$. It is not developed further here. Then, we essentially show
that typical paths (i.e. those having a non negligible contribution to the expectation) have no oriented edge read more than twice. This ensures that the expectation (\ref{largetrace}) only
depends on the variance of the entries $X_{ij}, i=1, \ldots, N, j=1, \ldots,p$. Universality of the expectation then follows.

\subsection{Number of self-intersections and odd marked instants in typical paths\label{subsec: nbrselfint}}We first give a technical Proposition which bounds the number of self intersections and give the approximate number of odd marked instants in typical paths.
In the following, we denote by $Z(k)$ the contribution of
paths with $k$ odd marked instants. We also denote the number of self-intersections on each line as
$M_1=\sum_{i\geq 2}(i-1)n_i\text{ and }M_2= \sum_{i\geq 2}(i-1)p_i.$ 

\bp \label{Prop: bound for self intersection}
There exists a positive constant $d_1$ such that the contribution of paths for which $M_1+M_2\geq d_1 \sqrt{ s_N}$
is negligible in the large-$N$-limit, whatever $1\leq k\leq s_N$ is.\\
And for any $\alpha, \alpha'$ such that $0<\alpha'<\frac{\sqrt \gamma}{1+\sqrt \gamma}<\alpha<1$, one has that
$$\sum_{k\leq \alpha's_N}Z(k) +\sum_{k\geq \alpha s_N}Z(k)=o(1) u_+^{s_N}.$$
\ep

\paragraph{Proof of Proposition \ref{Prop: bound for self intersection}}
We first give the proof of the first point of Proposition \ref{Prop: bound for self intersection}.
Denote by $Z(k,(\tilde n, \tilde p))$ the contribution of paths with $k$ odd marked instants and of type
$(\tilde n, \tilde p).$ Using (\ref{contrigal}), Remark \ref{rem : u+}
(and exactly the same arguments as in \cite{Sos} p. 34), one can see that, for $d_1$ large enough,
$$\sum_{k=1}^{s_N}\:\:\sum_{(\tilde n, \tilde p)/ \sum_{i\geq 2}(i-1)(n_i+p_i)\geq d_1 \sqrt{s_N}}Z(k,(\tilde n, \tilde p))=o(1) u_+^{s_N}.$$
\indent We now turn to the second statement. Let then $\alpha$ and $\alpha'$ be chosen as in
Proposition \ref {Prop: bound for self intersection}.
We assume that $N$ is large enough so that $\alpha'<\frac{\sqrt{\gamma_N}}{1+\sqrt{\gamma_N}}<\alpha.$
We now show that $\sum_{k\geq \alpha s_N}Z(k)<<u_+^{s_N}.$ Given any integer $k\leq s_N$, and using (\ref{contrigal}), one can show that there exists a constant $C_8>0$ independent
of $N$ and $k$ such that $Z(k)\leq N \gamma_N ^{k}\exp{\{C_8 N^{1/3}\}}\sigma^{2s_N}\mathbf{N}(s_N,k)$. Thus
\begin{eqnarray*}
&&\sum_{k \geq \alpha s_N } Z(k)\leq \sum_{k \geq \alpha s_N }\sigma^{2s_N}
 \mathbf{N}(s_N,k) \gamma_N^{k}N\exp{\{C_8 N^{1/3}\}}\crcr
 &&\leq \sigma^{2s_N}N C_{s_N}^{\hat k}C_{s_N}^{\hat k-1}\gamma_N^{\hat k}
\exp{\left\{C_8 N^{1/3}-C_7s_N\left (\alpha-\frac{\sqrt{\gamma_N}}{(1+\sqrt{\gamma_N})}\right )^2\right \}}\crcr
&& <<u_+^{s_N},
 \end{eqnarray*}
for $N $ large enough.
Similarly, the contribution of paths for which $k\leq\alpha' s_N$ is negligible in the large-$N$-limit. $\square$

\subsection{Asymptotics of $\mathbb{E}[\text{Tr}M_N^{s_N}]$}
In this subsection, we refine the estimate (\ref{contrigal}) and in particular deal with vertices of type $2$. Indeed, when summing (\ref{contrigal}) over $n_i, i\leq s_N$ and $p_i, i\leq s_N$, one can note that terms associated to vertices of type $2$  make the summation go to infinity.
To this aim, we shall control the number of vertices for which there is an ambiguity to continue the path at an unmarked instant. We shall also control the number of such vertices associated to edges passed four times or more.
Finally, we shall also show that amongst vertices of type
$3$, none belongs to edges passed more than twice, while there are no more complex self-intersections in typical paths.

\paragraph{}
From now on, given $\alpha' s_N\leq k\leq \alpha s_N$, we consider a path $P_k$ of type $(\tilde n, \tilde p)$ with $M_1=\sum_{i\geq 2}(i-1)n_i\leq d_1\sqrt{s_N}$ (resp. $M_2=\sum_{i\geq 2}(i-1)p_i\leq d_1\sqrt{s_N}$) self-intersections on the bottom (resp. top line). Our counting strategy is refined as follows.
Knowing the moments of self-intersection, we first choose the vertices occuring at the remaining marked moments.
One fills in the blanks of the path until the first moment of self-intersection is encountered. Then, 
one chooses the vertex which is repeated among the preceding ones (and repeat it if needed at the moment of second self-intersection and so on). We then proceed in the same way for subsequent vertices.\\
Assume that the instants of self-intersections have been chosen on each line:
$t_{jb,1}<t_{jb,2}<\cdots< t_{jb,n_2}$ for vertices of type 2 on the bottom line,
$t_{ju,1}<t_{ju,2}<\cdots< t_{ju,p_2}$ for vertices of type 2 on the top line,
$t_{jb,1}^{3,1}<t_{jb,2}^{3,1}<\cdots< t_{jb,n_3}^{3,1}$ for the first repetition of a
vertex of type 3 on the bottom line, $t_{jb,1}^{3,2}>t_{jb,1}^{3,1},
t_{jb,2}^{3,2}>t_{jb,2}^{3,1}\ldots$ for the second repetition of a vertex of type 3 on the
bottom line. 
We do not go deeper in the list of instants since it is exactly the same as in \cite{Sos} p. 724,
except that we make a distinction between instants marked on the top or bottom line (even if it is the same vertex).
The number of pairwise distinct vertices in the order of appearance on each line (top or bottom) occuring in the path is at most
\be N\prod_{i=1}^{s_N-k-M_1}(N-i) \prod_{i=1}^{k-M_2}(p+1-i). \label{nombresommets"distincts"} \ee 
Note that $(\ref{nombresommets"distincts"})\sim N^{s_N-M_1-M_2+1}\gamma_N^{k-M_2}
\exp{\{-\frac{(s_N-k)^2}{2N}-\frac{k^2}{2p}\}}.$\\
First, we focus on vertices of type 2. In the general case, there are ${ju,i}-i$ choices for the vertex occuring at the instant $t_{ju,i}$, since one chooses vertices occuring twice as marked instants. Assume first that the path is such that there
are no choices for closing any edge from such vertices at unmarked instants and that none belongs
to edges passed four times or more. Then choosing the instants and vertices of type 2 gives a contribution
at most of order
\begin{eqnarray*}
&&\sum_{1\leq t_{ju,1}<t_{ju,2}<\cdots< t_{ju,p_2}\leq s_N-k} \prod_{i=1}^{p_2}(j_{u,i}-i) \sum_{1\leq t_{jb,1}<t_{jb,2}<\cdots< t_{jb,n_2}\leq k}
\prod_{i=1}^{n_2}(j_{b,i}-i)\crcr
&&\leq \frac{1}{n_2!}\left ( \frac{(s_N-k)^2}{2}\right)^{n_2}\frac{1}{p_2!}\left ( \frac{k^2}{2}\right)^{p_2}.
\end{eqnarray*}
Such an estimate combined with formula (\ref{nombresommets"distincts"}) and Remark \ref{rem : u+} then ensures that the contribution of such paths to $\mathbb{E}[\text{Tr}(M_N/u_+)^{s_N}]$ is bounded.

\paragraph{}
Yet amongst vertices of type $2$ on the bottom or on the top line, at an unmarked instant where one closes an
edge with such a vertex as left endpoint, there might be a choice for closing the edge. Note that there are at
most three choices. An example of such a vertex is the distinguished vertex $4$ on Figure \ref{fig: cheminpk}, as the distinguished edge could have been $(4,4)$. Indeed, the two up edges with $4$ as marked vertex on the top line are read before the first time a down edge is closed starting from $4$. 
This leads to the notion of non-closed vertex. 

\bdefi A vertex $v$ of type 2 is said to be non-MP-closed if it is an odd (resp. even) marked instant and if there are more than one choice for closing an edge at an unmarked instant starting from this vertex on the top (resp. bottom) line.
\edefi

\brem The definition of non-MP-closed vertices differs from that of non-closed vertices in \cite{Sos}, essentially due to the distinction which is made between the top and bottom lines. 
\erem

Let $t$ be a given instant. Assume that the marked vertices before $t$ have been chosen and that, at the instant $t$, there is a non-MP-closed vertex. Then, by the definition of the trajectory and of non-MP-closed vertices, there are at most $x(t)$ possible choices for this vertex. This can be checked as in \cite{Si-So2}, p 122. In Lemma \ref{Lem: max x(t)} below, we show that $\max_{t} x(t) \sim \sqrt{s_N}$ in typical paths.

\paragraph{}
Apart from non-MP-closed vertices, a vertex of type $2$ can also belong to an edge that is read
four times in the path. To consider such vertices, we need to introduce other characteristics of
the path. Let $\nu_N( P)$ be the maximal number of vertices that can be visited at
marked instants from a given vertex of the path $P.$ Let also $T_N(P)$ be the maximal type of a
vertex in $P.$ Then, if at the instant $t$, one reads for the second time an
oriented up edge $e$, there are at most $2(\nu_N( P)+T_N(P))$ choices for the vertex occuring
at the instant $t$. Indeed, one shall look among the oriented edges already encountered in the path and of which one endpoint is the vertex occuring at the instant $t-1$ (see the Appendix in \cite{Si-So2} and \cite{PecheFeral} Section 5.1.2 e.g.).
It is an easy fact that paths for which $T_N(P)\geq A N^{1/3}(\ln N)^{-1}$
lead to a negligible contribution, if $A$ is large enough (independently of $k$). We prove at the end of this subsection, using Lemma \ref{Lem : nuN} stated below, that there exists $\epsilon>0$ small enough such that, for typical paths,
$$\nu_N( P)\leq s_N^{1/2-\epsilon} \text{ for any $\alpha's_N\leq k\leq \alpha s_N$}.$$
For vertices of type $i>2$, once the $i-1$ moments of self-intersection are fixed, one chooses at the first moment of self-intersection the vertex to be repeated amongst those already occurred in the path.\\
Assuming the above estimates on $\max x(t)$ and $\nu_N$ hold, we consider paths $P_k$ of type $(\tilde n, \tilde p)$ with $M_1:=\sum_{i\geq 2}(i-1)n_i\leq d_1\sqrt{s_N}$ and $M_2:=\sum_{i\geq 2}(i-1)p_i\leq d_1\sqrt{s_N}$ self-intersections respectively on the bottom or on the top line, $r_i$ non-MP-closed vertices of type $2$ ($i=1,2$) on the bottom and top lines and $q_i$ ($i=1,2$) vertices of type $2$ on the bottom or top line, visited at the second marked instant along an oriented edge already seen in the path. By Proposition \ref{Prop :closingandweight}, (\ref{nombresommets"distincts"}) and the above, their contribution to $\mathbb{E}[\text{Tr}M_N^{s_N}]$ is then at most of order (see also \cite{Sos}, p725)
\begin{eqnarray}
&&C\sigma^{2s_N}\mathbf{N}(s_N, k)N \gamma_N^{k}e^{\{-\frac{(s_N-k)^2}{2N}-\frac{k^2}{2p}\}}\: \mathbb{E}_k\Big [
\crcr &&\frac{1}{(n_2-r_1-q_1)!}\left(
\frac{(s_N-k)^2}{2N}\right)^{n_2-r_1-q_1}\frac{1}{r_1!}\left ( \frac{4(s_N-k)\max
\:x(t)}{N} \right)^{r_1}\crcr && \frac{1}{q_1!}\left (\frac{D_3(s_N-k)(\nu_N+T_N)}{N} \right)^{q_1}
\prod_{i\geq 3}\frac{1}{n_i!}\left (\frac{{C}^i(s_N-k)^i}{N^{i-1}}\right)^{n_i}\crcr &&
\!\!\frac{1}{(p_2-r_2-q_2)!}\left( \frac{k^2}{2p}\right)^{p_2-r_2-q_2}\frac{1}{r_2!}\left (
\frac{4k\max \:x(t)}{p} \right)^{r_2}\crcr && \frac{1}{q_2!}\left (\frac{D_4k(\nu_N+T_N) }{p}\right)^{q_2}
\prod_{i\geq 3}\frac{1}{p_i!}\left (\frac{{C}^ik^i}{p^{i-1}}\right)^{p_i}\Big ]
 \label{estcontrnunmax}
\end{eqnarray}
Here $D_3, D_4$ are positive constants independent of $k, p, N,$ and $s_N$, and $\mathbb{E}_k$ denotes the expectation with respect to the uniform distribution on $\mathcal{X}_{s_N,k}$. Here we have used the fact that 
$\sum_{x \in \mathcal{X}_{s_N,k}}f(x)=\mathbf{N}(s_N,k) \mathbb{E}_k (f(x))$ for any function $f\geq 0$.
Before considering paths in complete generality, we first restrict to paths with less than $d_1 \sqrt{s_N}$ self-intersections and no self-intersection of type strictly greater than $3$: $\sum_{i\geq 4}p_i+n_i=0.$\\
Let $Z_3(k)$ denote the total contribution of paths with $k$ odd marked instants such that $q:=q_1+q_2=0$, with no oriented edges read more than twice and satisfying the above conditions.

\bp \label{Prop:z3}There exists a constant $B_1>0$ independent of $N$ such that $Z_3:=\sum_{k=\alpha' s_N}^{\alpha s_N}Z_3(k)\leq B_1u_+^{s_N}.$
\ep

\paragraph{Proof of Proposition \ref{Prop:z3}}
From (\ref{estcontrnunmax}), we deduce that there exists a constant $D_o>0$ independent of $N,p,$ $k$ and $s_N,$ such that
\be \label{majoz3k} Z_3(k)\leq \sigma^{2s_N}\mathbf{N}(s_N, k)N\gamma_N^{k} \mathbb{E}_k\left ( \exp{\{6 \frac{\max x(t)s_N}{N}\}}\right) \exp{\{D_o \frac{s_N^3}{N^{2}}\}}.\ee
In Lemma \ref{Lem: max x(t)} proved below, we show that, given $a>0$, there exists $b>0$, independent of $N$, such that $\mathbb{E}_k\left ( e^{\{\frac{a \max x(t)}{\sqrt{s_N}}\}}\right)\leq b,$
$\forall \alpha's_N\leq k\leq \alpha s_N$.  This yields that $(\ref{majoz3k})\leq D_5\sigma^{2s_N}\mathbf{N}(s_N, k)N\gamma_N^{k},$
for some constant $D_5$ independent of $k$. Remark \ref{rem : u+} ensures that $Z_3:=\sum_{k}Z_3(k)=O(u_+^{s_N}).$ $\square$

\paragraph{}Assuming that there are no self-intersections of type greater than $3$, we can then show that paths for which $q=q_1+q_2\geq 1$ give a contribution of order $u_+^{s_N} \nu_N/\sqrt{ s_N}$ and thus there are no edges read more than twice (associated to vertices of type 2). We then proceed in the same way to show that there are no more than $\ln \ln N$ vertices of type $3$ in typical paths and that there are no oriented edges
read more than twice associated to vertices of type $3$.
It is then easy to deduce from the above result that paths with self-intersections of type $4$ or greater, or a marked origin, lead to a contribution of order $u_+^{s_N}s_N/N =o(1) u_+^{s_N}$. The detail is skipped.

\paragraph{}
Finally we investigate the total contribution of paths for which $\nu_N \geq s_N^{1/2-\epsilon}$ where $\epsilon >0$ is fixed (small). Denote by $Z_4$ such a contribution. We only indicate the tools needed to prove that 
$$Z_4:=\sum_{k=\alpha' s_N}^{\alpha s_N}Z_4(k)=o(1)u_+^{s_N},$$
since the detail of the proof is a rewriting of the arguments of the proof of Lemma 7.8 in \cite{PecheFeral}. 
To consider such paths, we introduce the following characteristic of the path, namely
$N_o:=r_1+r_2+\sum_{i\geq 3}in_i+ip_i.$
Assume then that $k, N_o,$ $q_1$ and $q_2$ are given.
We can then divide the interval $[0,2s_N]$ into $N_o$ sub-intervals, so that inside an interval, there are only closed vertices of type 2 or no self-intersection. Then there is no choice for closing the edges inside these sub-intervals.
Assume that a vertex $v$ is the starting point of $\nu_N$ up edges.  Then, there is a time interval $[s_1, s_2]$
during which the trajectory of $P_k$ comes $\nu'_N:=\frac{\nu_N}{2N_o}$ times to the level $x_o$ (of $v$) and never goes below. Denote by $\Gamma (\nu'_N)$ the event that there exists such an interval in a trajectory and let $\mathbb{P}_k$ denote the uniform distribution on $\mathcal{X}_{s_N,k}.$ In Lemma \ref{Lem : nuN}, we show that the probability of such an event decreases as %
$$\max_{1\leq k\leq s_N}\mathbb{P}_k \left (\Gamma (\nu'_N)\right)\leq A_1s_N^2\exp{\{-A_2\nu'_N\}}, $$ 
for some positive constants $A_1, A_2$ independent of $s_N$. Using Lemma \ref{Lem: max x(t)}, one can also show that there exists a constant $A>0$ such that $\max x(t)\leq A N^{1/6}\sqrt{s_N}$ in any non-negligible path. Using these estimates, formula $(\ref{estcontrnunmax})$, and Lemma \ref{Lem : nuN} proved below, one can then copy the arguments used
in \cite{PecheFeral} Lemma 7.8, to deduce that
$Z_4=o(1)u_+^{s_N}.$\\ 
This finishes the proof that typical paths have a non-marked origin, vertices of type $3$ at most (and less than $\ln \ln N$ of type $3$), less than $d_1\sqrt{s_N}$ self-intersections and no edges read more than twice.
The proof of Theorem \ref{theo: unimoments} is
completed once Lemmas \ref{Lem: max x(t)} and \ref{Lem : nuN} are proved. $\square$

\subsection{Technical Lemmas \label{Subsec: techlemmas}}
In this subsection, we prove the results used in the previous subsection on characteristics
of typical paths. The first quantity of interest here is the maximum level reached by the
trajectory of a path, namely $\max_{t \in [0, 2s_N]} x(t)$. We shall show that it roughly behaves as
$\sqrt{s_N}$ in typical paths. The second one is the maximal number of vertices visited from
a given vertex, $\nu_N( P),$ which should not grow faster than $s_N^{\delta}$ for any
power $\delta$ (but we get a weaker bound).

\paragraph{}
We shall now prove the announced estimate for  $\max_{t\in[0, 2s_N]} x(t)$, where $x(t)$ denotes the level
of the trajectory $X_k$ associated to a path $P=P_k$. Let $a$ be some constant independent of $N$ and $k$
and denote by $\mathbb{E}_{k}$ the expectation with respect to the uniform distribution $\mathbb{P}_k$ on the set of trajectories
$\mathcal{ X}_{s_N,k}$ of length $2s_N$ with $k$ up edges at odd instants.

\bl \label{Lem: max x(t)} Given $a >0,$ there exists $b>0$ independent of $s_N$ (and $N$) such that
$$\max_{\alpha' s_N\leq k\leq \alpha s_N}\mathbb{E}_{k}\exp{\left \{\frac{a\max x(t)}{\sqrt{ s_N}}\right \}}\leq b.$$\el

\begin{figure}[htbp]
 \begin{center}
 \begin{tabular}{c}
 \epsfig{figure=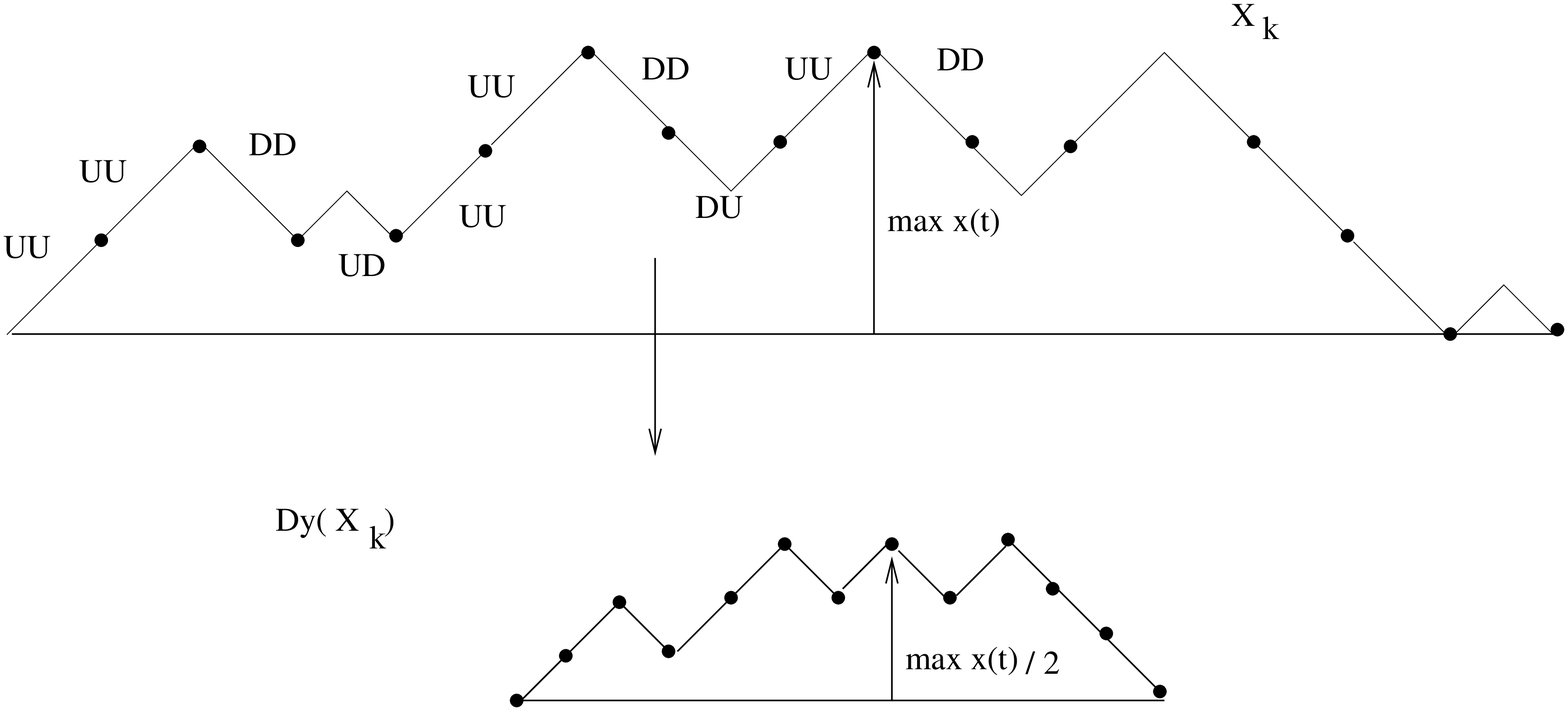,height=6cm,width=10cm, angle=0}
 \end{tabular}
 \caption{A trajectory $X_k$ and the associated trajectory $Dy(X_k)$.
 \label{fig: UUDD}}
 \end{center}
\end{figure}

\paragraph{Proof of Lemma \ref{Lem: max x(t)}}
It was proved in \cite{Si-So2} that the above result holds if one replaces $\mathbb{E}_{k}$ with the expectation with respect to the uniform distribution on Dyck paths (no constraint on $k$) of length $2s_N.$ We will call on this result to prove Lemma \ref{Lem: max x(t)}. To this aim, we cut the Dyck path $X_k$ into $2-$steps, so that there are 4 types of basic 2-steps : $UU$, $UD$, $DD$ and $DU$ ($D$ stands here for down, $U$ for up). It is an
easy fact that the number of $UU$ steps equals that of $DD$ steps. Let then $l$ be the number of
$UU$ steps (and $DD$ steps), $k_2$ be those of $DU$ and $k_3$ be those of $UD$ steps. Then,
\be \label{estiml}
2l+k_3+k_2=s_N, \: l+k_2=s_N-k,\: l+k_3=k. 
\ee
As a step $UD$ or $DU$ brings the path to
the same level, it is easy to see that the steps $UU$ and $DD$ are arranged in such a way that
they form a Dyck path (if we identify a $UU$ step with an up step and a $DD$ step with a down step) of length $2l$. We denote $Dy(X_k)$ this sub-Dyck path associated to the trajectory $X_k$ (see Fig. \ref{fig: UUDD}).

We now explain how to build a trajectory $X_k$ given $Dy(X_k)$ of length $2l$ and $s_N-k-l$ (resp. $k-l$) $DU$ (resp. $UD$) steps.
To construct $X_k$ from $Dy(X_k)$, one has to insert ``horizontal'' steps, namely $DU$ and $UD$ steps, in a particular way. Note that two distinct insertions lead to two different trajectories.
The sole constraint is to insert steps $DU$ when the path $Dy(X_k)$ is at a level
greater than or equal to one. This is the reason why we enumerate Dyck paths with $2l$ steps
according to the number of times they come back to the level $0$. Call $\sharp Dyck(l,Q)$ the
number of Dyck paths with $2l$ steps and $Q$ returns to $0$. We then have to
insert $s-k-l$ horizontal $DU$ steps into $2l-Q$ boxes. Then we can insert the $UD$ steps
arbitrarily. This yields that
\begin{eqnarray}
&&\mathbf{N}(s_N,k)=\sum_{l=0}^{k\wedge (s_N-k)}\:  \sum_{Q=0}^l \: \sharp
Dyck(l,Q)C_{l-Q+s-k-1}^{s-k-l}C_{s}^{k-l}.\label{dypathretours0}
\end{eqnarray}
From this construction, it is easy to see that the maximum level reached by the trajectory $X_k $ is twice the maximum $(+1)$
reached by the sub-path $Dy(X_k).$
Let then
$\mathcal{Y}_{l,Q}$ denote the set of Dyck paths of length $2l$ with $Q$ returns to $0.$ 
We do not consider the degenerate case where $Dy(X_k)=\emptyset$, which corresponds to the trajectory obtained with $UD$ steps only. Then one has that
\begin{eqnarray}
&& \mathbb{P}_k(\max X_k (t)=r )\crcr %
&&\!\!\!\!\!\leq  \sum_{l=1}^{s_N/2}\sum _{Q=0}^{s_N} \mathbb{P}_k \left(\max Dy(X_k) =r/2 \big | Dy(X_k)\in
\mathcal{Y}_{l,Q}\right ) \mathbb{P}_k\left (Dy(X_k)\in \mathcal{Y}_{l,Q}\right ) +\crcr
&&\!\!\!\!\! \sum_{l=1}^{s_N/2}\sum _{Q=0}^{s_N} \mathbb{P}_k \left(\max Dy(X_k) =\frac{r-1}{2} \big |
Dy(X_k)\in \mathcal{Y}_{l,Q}\right ) \mathbb{P}_k\left (Dy(X_k)\in \mathcal{Y}_{l,Q}\right )
.\crcr
&&\label{estimeemaxenq}
\end{eqnarray}
Let $\mathbb{P}_{l,Q}$  denote the uniform distribution on $\mathcal{Y}_{l,Q}$.
Let also $a_o>0$ (small) be given. It can easily be
inferred from \cite{Si-So2}, p. 11 (see also \cite{Sos} and \cite{PecheFeral}, Lemma 7.10) that there exist positive constants $a_1$, $a_2,$ independent of $l$ and $Q$ such that, if $r'\geq a_o \sqrt{l}$, one has that 
\be \label{estxsos}\mathbb{P}_{l,Q} (\max x (t)= r')\leq \frac{a_1}{\sqrt l}\exp{\{-\frac{a_2{r'}^2}{l}\}}.\ee
Thus, inserting (\ref{estxsos}) in (\ref{estimeemaxenq}), we deduce that there exist some positive constants
$a_3, a_4$ independent of $k,$ $s_N$ and $N$ such that, provided $r \geq a_o \sqrt{s_N}$ and for any $\alpha' s_N\leq k\leq \alpha s_N,$ $\mathbb{P}_k(\max X_k (t)=r )\leq \frac{a_3}{\sqrt{s_N}}\exp{\{-\frac{a_4{r}^2}{s_N}\}}.$ This yields Lemma
\ref{Lem: max x(t)}. $\square$

\paragraph{}The second estimate we need is a suitable bound on $\nu_N( P_k)$. Recall that $\Gamma (\nu'_N)$ denotes the event that the trajectory of a path comes back from above $\nu'_N$ times to some level $x_o$. 

 \bl \label{Lem : nuN}%
There exist positive constants $A_1, A_2$ independent of $k, N, p$ such that%
\be \label{estinitnu}\max_{1\leq k\leq s_N}\mathbb{P}_k \left (\Gamma (\nu'_N)\right)\leq
A_1s_N^2\exp{\{-A_2\nu'_N\}}. \ee
  \el

\paragraph{Proof of Lemma \ref{Lem : nuN}}
Let $[s_1, s_2]$ be an interval such that $x(t_1)=x(t_2)=x_o$ for some $x_o\geq 0$
and $x(t)\geq x_o, \forall t \in[s_1,s_2]$ and for which there exists $ s_1<t_1<t_2<\cdots<t_{\nu'_N}\leq s_2$ such that $x(t_i)=x_o$. 
We first consider the case where $s_1$ and $s_2$ are even instants (then $x_o$ is also even).
Modifications to be made in the case where they are odd will be indicated at the end of the proof.
The instants $t_i$ are then called instants of returns from above to $x_o$ of the trajectory $X_k$.
Set now $Y_o$ to be the Dyck path of length $s_2-s_1$ defined by $y_o(t)=x(t+s_1)-x_o, t\in [0, s_2-s_1]$.
Then the returns from above to $x_o$ correspond to returns to $0$ of $Y_o$.
Now, the returns to $0$ of $Y_o$ can either be made using $UD$ steps or correspond to a return of the sub-trajectory $Dy(Y_o)$ to this level. Thus, either the number of $UD$ steps is large or the number of returns of $Dy(Y_o)$ to level $0$ is large. We shall show that in both cases, (\ref{estinitnu}) holds. Thanks to Proposition \ref{Prop: bound for self intersection}, it is enough to consider trajectories $X_k$ for which $\alpha's_N\leq k\leq \alpha s_N.$ The proof of Lemma \ref{Lem : nuN} is divided into three steps.

\paragraph{Step 1}
We first show that there exist positive constants $C'_7, C'_8$ independent of $N$ and $k$ such that, provided $\alpha's_N\leq k \leq \alpha s_N,$
\be \label{eststep1}\mathbb{P}_k\left ( \text{ $X_k$ has $\eta_N$ consecutive $UD$ steps }\right)\leq C'_7s_N^2 \exp{\{-C'_8 \eta_N\}}.\ee
Assume that there exists a time interval $[s'_1,s'_2]$ with $\eta_N$ consecutive $UD$ steps only.
Given even instants $s'_1$ and $s'_2$ (with $s'_2-s'_1=2\eta_N$),
the proportion of trajectories $X_k$ that have $\eta_N$ steps $UD$ in $[s'_1, s'_2]$ is at most
\begin{equation}
\frac{\frac{1}{s_N-\eta_N}\left (C_{s_N-\eta_N}^{k-\eta_N}\right)^2}{\frac{1}{s_N}\left (C_{s_N}^{k}\right)^2}
\leq C_9\left ( \frac{k(k-1)\cdots (k-\eta_N+1)}{s_N(s_N-1)\cdots (s_N-\eta_N+1)}\right)^2
\!\leq C_9 \alpha^{2\eta_N},\label{estimeenu'ninit}
\end{equation}
for some constant $C_9>0.$ This readily yields (\ref{eststep1}).

\paragraph{Step 2}We consider the case where the number of returns to 
level $0$ made by the associated path $Dy(Y_o)$ is large.
It was proved in \cite{Si-So2} that, if $\mathbb{P}_{s}$ denotes the uniform distribution on the set of
Dyck paths $Y$ with length $2s$, then there exist constants $C_{10}, C_{11}$ independent of $s$ such that 
\begin{eqnarray} \label{estimeenu'n2}
&&\mathbb{P}_s(\exists s'_1,s'_2: Y \text{ has $\eta_N$ returns from above to the level
 $x_o$ in $[s'_1,s'_2]$})
\crcr
&&\leq C_{10}s^2\exp{\{-C_{11}\eta_N\}}.
\end{eqnarray}
Denote by $Q$ the number of sub-Dyck paths $\tilde Y_i, i=1, \ldots Q,$ of $Dy(Y_o)$ starting and ending at level $0.$ 
From the above result and (\ref{dypathretours0}), one can deduce that there exists constants $C_{12}, C_{13}>0$, independent of $s_N$ and $k$, such that 
\be
P(Q=\eta_N)\leq C_{12}s_N^2 \exp{\{-C_{13}\eta_N\}}. \label{estimeenu'_ninit2}
\ee

\paragraph{Step 3}We can now turn to the proof of Lemma \ref{Lem : nuN}.
A trajectory $X_k$ coming back $\nu'_N $ times to the level $x_o$ during $[s_1, s_2]$ can be described as
follows. Denote by $Q$ the number of sub-Dyck paths $Y_i, i=1, \ldots Q,$ going from level $x_o$
to $x_o$ and starting with a $UU$ step and ending with a $DD$ step. Denote by $l_i,i=1, \ldots Q,$ the respective length of these sub-Dyck paths.
Then these sub-Dyck paths are interspaced by $\nu'_N-Q$ $UD $ steps that split in at most $Q+1$ sequences.
Let $\nu_N^i$ ($i\leq Q+1$) be the respective lengths of these disjoint sequences of $UD$ steps from $x_o$ to $x_o$.
Using the estimates of Step 1 and Step 2, there exist constants $C'_{12}, C'_{13}>0$ such that for any constants $A, A'>0$ (fixed later),
\begin{eqnarray}
&&\mathbb{P}_k\left ( \exists \:\nu_N^i \geq \nu'_N/A \right)\leq C'_{12}s_N^2\exp{\{-C'_{13}\nu'_N/A\}},\: \forall \,1\leq k\leq s_N,\crcr
&&\mathbb{P}_k\left ( Q\geq \nu'_N/A'\right ) \leq C'_{12}s_N^2\exp{\{-C'_{13}\nu'_N/A'\}},\: \forall \,1\leq k\leq s_N.
\end{eqnarray}
Thus it is enough to consider trajectories such that $A\leq Q\leq \nu'_N/A'.$ To count these
trajectories, we study their structure in more detail. Set $L_o=(s_2-s_1)/2$ and let then $k_o<k$
be the number of odd marked instants of the sub-trajectory inside the interval $[s_1, s_2]$. The
remaining trajectory $x(t)$, $t \in [0, 2s_N]\setminus [s_1, s_2],$ is then a Dyck path of length $2s_N-2L_o$ with $k-k_o$ odd up steps. Given $s_1$ and $s_2,$ and in order to count the number of such trajectories $X_k,$ we first re-order the paths $Y_i$ and $UD$ steps inside the interval $[s_1, s_2]$ as follows. We first read the Dyck paths $Y_i, i=1, \ldots Q,$ and then read all the $UD$ steps.\\
Fix some $0<\epsilon <(1-\alpha)/2$ (small). Assume first that 
\be L_o-k_o<(1-\epsilon)s_N+\epsilon (L_o-(\nu'_N-Q))-k. \label{condko}\ee
The latter condition ensures, as $k\leq \alpha s_N$ for some $\alpha<1$, that $k-k_o+\nu'_N-Q\leq (1-\epsilon)(s_N-L_o+\nu'_N-Q)$. Thus, we can apply Step 1 to the sub-trajectory obtained from $x$ by erasing the sub-paths $Y_i, i=1, \ldots, Q$.
Then, given $Q, k_o$, $s_1$ and $s_2$, the number of trajectories of length $L=2s_N-2L_o +2(\nu'_N-Q)$ with $k-k_o+\nu'_N-Q$ odd up steps and that have $\nu'_N-Q$ $UD$ steps
between $[L-(s_N-s_2)-2(\nu'_N-Q), L-(s_N-s_2)]$ is of order
$$ C_{14} \exp{\{-C_o' (\nu'_N-Q)\}}\times N_2,\text{ if }N_2=\sharp\{\mathcal{X}_{ 2s_n-2L_o +2(\nu'_N-Q),\, k-k_o+\nu'_N-Q}\}.$$
Here $C_{14},C_o'$ are some positive constants independent of $s_N, k$ and $N.$ Note that the constant $C'_o$ depends only on $\epsilon$ and $\alpha$. Then, the number of
Dyck paths of length $\sum_{i=1}^Q l_i=2L_o-2(\nu'_N-Q)$, with $k_o-(\nu'_N-Q)$ odd up steps and coming back $Q$ times to the
level 0 using $DD$ steps is at most of order
$$ C_{15} \exp{\{-C_o'' Q\}}\times N_1,
\text{ if }N_1=\sharp\{ \mathcal{X}_{\sum_{i=1}^Q l_i ,\, k_o-(\nu'_N-Q)}\}.$$ As above $C_{15}, C_o''$
are positive constants independent of $N, s_N$ and $k.$\\
Finally the number of ways to order the paths
$Y_i$ and the $UD$ steps inside the interval $[s_1,s_2]$ is equal to the number of
ways to write $\nu'_N-Q $ as a sum of $Q+1$ integers. There are $C_{\nu'_N}^{Q}$ such ways.\\
Thus the number of trajectories $X_k$ coming $\nu'_N$ times to some level $x_o$ never falling
below is at most, if $C_o=\min \{C_o', C_o''\},$
\begin{eqnarray}
&&\sum_{0\leq s_1<s_2\leq s_N }\sum_{Q=A}^{\nu'_N/A'} \sum_{k_o\leq k}\: C_{\nu'_N}^{Q}C_{16} \exp{\{-C_o\nu'_N \}} N_1 N_2\crcr &&\leq C_{16}s_N^2
\sum_{Q=A}^{\nu'_N/A'}C_{\nu'_N}^{Q}\exp{\{-C_o \nu'_N\}}\mathbf{N}(s_N,k),\label{sommenunQ}
\end{eqnarray}
since, $L_o, Q,k$ and $\nu'_N$ being fixed, $\sum_{k_o}N_1 N_2 \leq \mathbf{N}(s_N,k)$. This yields the following estimate:
$$\mathbb{P}_k \left (X_k \text{ has }\nu'_N \text{ returns to }0, \:A\leq Q\leq \frac{\nu'_N}{A'}
\right ) \leq \nu'_N s_N^2e^{-C_o\nu'_N}C_{\nu'_N}^{\nu'_N/A'}.$$
We can then choose $A'$ large enough so that there exists a constant $C_{18}>0,$ independent of
$N, k$ and  $s_N,$ such that
$$\nu'_Ne^{-C_o\nu'_N}C_{\nu'_N}^{\nu'_N/A'}\leq C_{18}\exp{\{-C_o\nu'_N/2\}}.$$
This yields Lemma \ref{Lem : nuN} if (\ref{condko}) is satisfied.\\ 
Assume now that (\ref{condko}) is not satisfied. Then necessarily $k_o\leq (\alpha+\epsilon)L_o\leq(\alpha+1)L_o/2.$
Thus the number of trajectories $Y_o$ coming $\nu'_N$ times to some level $x_o$ with $Q$ returns made using $DD$ steps is at most 
$$C_{\nu'_N}^{Q} \mathbf{\tilde N}(L_o-(\nu'_N-Q),k_o-(\nu'_N-Q) , Q),$$
where $ \mathbf{\tilde N}(L_o-(\nu'_N-Q),k_o-(\nu'_N-Q) , Q)$ is the number of Dyck paths of length $2L_o-2(\nu'_N-Q)$, with
$Q$ returns to $0$ made only with $DD$ steps and admitting $k_o-(\nu'_N-Q)$ odd up steps.
From Step 2, one has that there exists a constant $C_{20}$ (independent of $k_o,L_o$) such that
$$\mathbf{\tilde N}(L_o-(\nu'_N-Q),k_o-(\nu'_N-Q) , Q)\leq e^{\{-C_{20}Q\}}\mathbf{N}(L_o-(\nu'_N-Q),k_o-(\nu'_N-Q) ).$$
As $k_o\leq (\alpha+\epsilon)L_o$, there also exists $C_{21}>0$ (depending on $\epsilon$ and $\alpha$ only) such that
$\mathbf{N}(L_o-(\nu'_N-Q),k_o-(\nu'_N-Q) )\leq \exp{\{-C_{21} (\nu'_N-Q)\}}\mathbf{N}(L_o, k_o).$
The end of the proof is as above. This finishes the proof of Lemma \ref{Lem : nuN} in the case where $s_1$ and $s_2$ are even.\\
To consider the case where $s_1$ and $s_2$ are odd, one can then use exactly the same arguments as above, up to the following modifications. An odd marked instant of $Y_o$ simply defines an even marked instant of $X_k$. Then it is an easy task to show that Step 1 holds if one replaces $\alpha$ with $1-\alpha'$ in (\ref{estimeenu'ninit}). Step 2 and Step 3 can then be obtained using arguments as above. $\square$

\paragraph{}
This finally completes the proof of Theorem \ref{theo: unimoments}.

\section{The case where $p, N \to \infty$ and $N/p\to 0$ \label{Sec: gammainf}}
In this section, we consider the sequence of random sample covariance matrices $M_p=\frac{1}{p}XX^*$ instead of $M_N$. We also consider moments of traces of $M_p$ to some powers 
$s_N\sim \sqrt{\gamma_N}N^{2/3}.$
The above scaling readily comes from Theorem \ref{theoNEK} proved by \cite{NEK}.
In particular, one has that \be\label{espmp}\mathbb{E}\Bigl(\text{Tr} M_p^{s_N}\Bigr)=\frac{1}{\gamma_N^{s_N}}\mathbb{E}\Bigl(\text{Tr} M_N^{s_N}\Bigr).\ee 
We here prove the following universality result.
Let $K$ be some given integer and $c_1, \ldots , c_K$ be positive constants chosen in a fixed compact interval $J$ of $\mathbb{R}_+^*$. 
Consider sequences $s_N^{(i)}, i=1, \ldots, K,$ such that $\lim_{N \to \infty}\frac{s_N^{(i)}}{\sqrt{\gamma_N}N^{2/3}}=c_i.$
\bt \label{theo: unigammainfini2}
Let $ v_+=\sigma^2\left (1+\frac{1}{\sqrt{\gamma_N}}\right )^2.$ There exists $C_1(K)>0$ such
that
\begin{eqnarray*}
&&\mathbb{E} \left  [ \prod_{i=1}^K \emph{Tr}  \left(\frac{XX^*}{pv_+} \right)^{s_N^{(i)}}\right]\!\!\leq C_1(K) \text{ and }\crcr
&&\mathbb{E}\left[ \prod_{i=1}^K \emph{Tr}  \left(\frac{XX^*}{pv_+}\right)^{s_N^{(i)}}
\right]=\mathbb{E}\left[ \prod_{i=1}^K \emph{Tr}  \left(\frac{X_GX_G^*}{pv_+}\right)^{s_N^{(i)}}\right](1+o(1)).
\end{eqnarray*} \et

\paragraph{}The proof of Theorem \ref{theo: unigammainfini2} is the object of the whole section. We only consider the case where $K=1$ and $s_N$ is some sequence such that
\be
\lim_{N \to \infty}\frac{s_N}{\sqrt{\gamma_N}N^{2/3}}=c,\text{ for some real $c>0$}.
\ee
As in the preceding section, we establish that the typical paths have no edges read more than twice. This ensures that the leading term in the asymptotic expansion of  $\mathbb{E} \left [\text{Tr} M_p^{s_N}\right ]$ is the same as that for Wishart ensembles. The idea of the proof is very similar to that of the preceding section, but requires some minor modifications. This is essentially due to the discrepancy between marked vertices on the bottom and top lines, due to the fact that $p>>N.$\\
In this section, $C$, $C_i, C'_i,D_i,B_i, i=0, \ldots, 9,$ denote some positive constants independent of $N, p,k$ and $s_N$ whose value may vary from line to line (and from the preceding sections). 

\subsection{Typical paths \label{subs: typicalpathsginf}}
We now state the counterpart of the second point in Proposition \ref{Prop: bound for self intersection} in the two following Propositions.
Let $Z_{\infty, 0}=Z_{\infty, 0}(A)$ be the sub-sum corresponding to the contribution to (\ref{espmp}) of paths for which $k\leq s_N(1-\frac{A}{\sqrt{\gamma_N}})$, for some $A>0$ to be fixed.

\bp \label{Prop: Zinfini0}There exists $A>0$ such that $Z_{\infty, 0}=o(1)v_+^{s_N}.$
\ep

\paragraph{Proof of Proposition \ref{Prop: Zinfini0}}
By (\ref{contrigal}), one has that
\begin{eqnarray}\label{estzinf}
&Z_{\infty, 0}\leq &\sum_{k=1}^{s_N (1-\frac{A}{\sqrt{\gamma_N}})}  \frac{N}{s_N}\frac{k}{s_N-k+1}
\left ( C_{s_N}^k\right)^2 \frac{2\sigma^{2s_N}}{\gamma_N^{s_N-k}}\crcr
&& \sum_{(\tilde n, \tilde p)}\prod_{i\geq 2}
\frac{1}{n_i!}\left( \frac{C^i (s_N-k)^i}{N^{i-1}}\right)^{n_i}\prod_{i\geq 2}\frac{1}{p_i!}
\left( \frac{C^i k^i}{p^{i-1}}\right)^{p_i}.\end{eqnarray}
Now there exists a constant $C_1>0$ such that, for any $1\leq k\leq s_N$,
\be \sum_{p_i, 1\leq i\leq k}\:\prod_{i\geq 2}\frac{1}{p_i!}\left( \frac{C^i s_N^i}{p^{i-1}}\right)^{p_i}
\leq \exp{\{C_1N^{1/3}\}} . \label{sumpgammainf}\ee
Similarly, using the fact that $\sum_{i=1}^{s_N-k}in_i=s_N-k$, we deduce that
\begin{eqnarray}
&&\frac{1}{\sqrt{\gamma_N}^{s_N-k}}\sum_{n_i, 1\leq i\leq s_N-k}\: \: \,\prod_{i\geq 2}\frac{1}{n_i!}
\left( \frac{C^i (s_N-k)^i}{N^{i-1}}\right)^{n_i}\crcr
&&\leq \frac{1}{\sqrt{\gamma_N}^{s_N-k}}\sum_{n_i, 1\leq i\leq s_N-k}\:\,\prod_{i\geq 2}\frac{1}{n_i!}
\left( \frac{C^i s_N^i}{N^{i-1}}\right)^{n_i}\left(\frac{s_N-k}{s_N} \right)^{in_i}\crcr
&&\leq \exp{\{C_2 N^{1/3}\}} \left ( \frac{s_N-k}{s_N}\right)^{s_N-k},\label{sumNgammainf}
\end{eqnarray}
as $s_N-k>A s_N \gamma_N^{-1/2}$ and provided $A>2$. 
Inserting (\ref{sumNgammainf}) and (\ref{sumpgammainf}) in (\ref{estzinf}) yields that
\begin{equation}
Z_{\infty, 0}\leq \sigma^{2s_N}\sum_{k=1}^{s_N (1-\frac{A}{\sqrt{\gamma_N}})} \frac{N}{s_N}\frac{k}{s_N-k+1}
\left ( C_{s_N}^k\right)^2 
\left( \frac{s_N-k}{\sqrt{\gamma_N}s_N}\right)^{s_N-k} e^{\{C_3 N^{1/3}\}}. \label{uppbdZinf2}
\end{equation}
We then deduce the following upper bound. For $N$ large enough, one has 
\begin{eqnarray*}
&\!\!\!\!\!\!\!(\ref{uppbdZinf2})&\leq C_4\sigma^{2s_N}N \sum_{k\leq s_N(1-\frac{A}{\sqrt{\gamma_N}})}
\frac{e^{\{C_3N^{1/3}\}}}{\sqrt{\gamma_N}^{s_N-k}}\left ( \frac{2s_N^2e^2}{(s_N-k)^2}\right)^{s_N-k}
\left ( \frac{s_N-k}{s_N}\right)^{s_N-k} \crcr
&&\leq C_4\sigma^{2s_N} \exp{\{C_4 N^{1/3}\}}\sum_{k\leq s_N(1-\frac{A}{\sqrt{\gamma_N}})} \left ( \frac{2e^2}{A}\right)^{s_N-k}=o(1) v_+^{s_N},
\end{eqnarray*}
where in the last line we have chosen $A>4e^2.$
This finishes the proof of Proposition \ref{Prop: Zinfini0}. $\square$

\paragraph{}
Given $0<\epsilon<1/2$, we also consider the contribution $Z_{\infty,1}=Z_{\infty,1}(\epsilon)$ of paths for which $k \geq s_N(1- \frac{\epsilon}{\sqrt{\gamma_N}})$.
\bp \label{Prop:zinfini1} There exists $0<\epsilon<1/2$
such that $Z_{\infty,1}=o(1)v_+^{s_N}.$ \ep

\paragraph{Proof of Proposition \ref{Prop:zinfini1}}
By (\ref{contrigal}) and as $s_N-k =O(N^{2/3})$, one has that
\begin{eqnarray}
&Z_{\infty,1}&\leq 2\sigma^{2s_N}\sum_{k=s_N(1- \frac{\epsilon}{\sqrt{\gamma_N}})}^{s_N}N \mathbf{N}(s_N,k)
\exp{\{(C_1+C'_2)N^{1/3}\}}\gamma_{N}^{-(s_N-k)}\crcr 
&&\leq \sigma^{2s_N}C'_4N \exp{\{C'_4N^{1/3}\}}
\mathbf{N}(s_N,[s_N-[\frac{\epsilon s_N}{\sqrt{\gamma_N}}+1])\gamma_{N}^{-[\frac{\epsilon s_N}{\sqrt{\gamma_N}}]-1}\crcr
 &&\leq C'_4N
\exp{\{C'_4N^{1/3}\}}\exp{\{-C_5N^{2/3}/8\}}v_+^{s_N},
\end{eqnarray}
provided $\epsilon <1/2.$
This is enough to ensure Proposition \ref{Prop:zinfini1}. $\square$

\paragraph{}Set now $I_N=[s_N(1-\frac{A}{\sqrt{\gamma_N}}), s_N(1-\frac{\epsilon}{\sqrt{\gamma_N}})]$.
Thanks to Proposition \ref{Prop: Zinfini0} and Proposition \ref{Prop:zinfini1}, typical paths are such that $k\in I_N$. This implies in particular that
$\frac{(s_N-k)^2}{N}=O(N^{1/3})$  and $\frac{k^2}{p}=O(N^{1/3}).$
Using the fact that $\sum_{k \in I_N}\sigma^{2s_N}\mathbf{N}(s_N, k)N \gamma_N^{k-s_N}=O(v_+^{s_N})$, it is easy to deduce from (\ref{contrigal}), that it is enough to consider paths for which 
$M_1+M_2=\sum_{i\geq 2}(i-1)(n_i+p_i)\leq d_1 N^{1/3}$ for some constant $d_1$ independent of $k\in I_N$, $N$ and $p$. 
For such paths, denote by $Z_{\infty}(k)$ the contribution of paths with $k$ odd marked instants.
We now prove the following Proposition yielding Theorem \ref{theo: unigammainfini2}.
\bp \label{Prop: derestgammainf}There exists $D_1\!>\!0$ such that 
$\sum_{k \in I_N}Z_{\infty}(k)\leq D_1v_+^{s_N}.$
Furthermore, the contribution of paths admitting either an edge read more than twice, or more than $\ln \ln N$ vertices of type $3$, or a vertex of type 4 or greater, or a marked origin is negligible in the large-$N$-limit.
\ep
\brem Considering paths $P_k, k \in I_N,$ admitting only vertices of type 2 at most and no non-MP-closed vertices, one can also deduce from the subsequent proof of Proposition \ref{Prop: derestgammainf}, that there exists $D_2>0$ such that $\sum_{k \in I_N}Z_{\infty}(k)\geq D_2 v_+^{s_N}.$
\erem

\paragraph{Proof of Proposition \ref{Prop: derestgammainf}}First, one can state the counterpart of Formula (\ref{estcontrnunmax}). Due to the different scales
$s_N-k =O(N^{2/3}),$ while $k=O(\sqrt{\gamma_N}N^{2/3})$, we need in this section to
 distinguish vertices being the left endpoint of an up edge according to the parity of the corresponding instant.
Define $\nu_{N,o}(P_k)$ (resp. $\nu_{N, e}(P_k)$ to be the maximum number of vertices
visited (at marked instants) from a vertex of the path occuring at odd instants (resp. even instants). Then,
\begin{eqnarray}
&Z_{\infty}(k)&\leq   C\sigma^{2s_N}\mathbf{N}(s_N, k)N \gamma_N^{k-s_N}e^{\{-\frac{(s_N-k)^2}{2N}-\frac{k^2}{2p}\}}\crcr
&& \sum_{n_2, r_1,q_1, n_3,\ldots, n_{s_N-k} }\:\:\sum_{p_2, r_2,q_2, p_3,\ldots, p_{k} }\mathbb{E}_k\Big [\frac{\left(
(s_N-k)^2/(2N)\right)^{n_2-r_1-q_1}}{(n_2-r_1-q_1)!}\crcr
&&\frac{1}{r_1!}\left ( \frac{4(s_N-k)\max
\:x(t)}{N} \right)^{r_1}\frac{1}{q_1!}\left (\frac{D_3(s_N-k)(\nu_{N,o}+T_N)}{N} \right)^{q_1} 
\crcr && 
\frac{\left( k^2/(2p)\right)^{p_2-r_2-q_2}}{(p_2-r_2-q_2)!}\frac{1}{r_2!}\left (
\frac{4k\max \:x(t)}{p} \right)^{r_2}\frac{1}{q_2!}\left (\frac{D_4k(\nu_{N,e}+T_N) }{p}\right)^{q_2} 
\crcr
&& \prod_{i\geq 3}\frac{1}{n_i!}\left (\frac{{C}^i(s_N-k)^i}{N^{i-1}}\right)^{n_i}\prod_{i\geq 3}\frac{1}{p_i!}\left (\frac{{C}^ik^i}{p^{i-1}}\right)^{p_i}\Big ].
\label{contrikinfini}
\end{eqnarray}
One still has that $T_N<A''N^{1/3}/\ln N$ for some $A''>0$ in typical paths (independently of $k\in I_N$).
This ensures that the analysis performed for the case where $\lim_{N \to \infty}\gamma_N<\infty$ can be copied,
 provided the counterparts of technical lemmas of Subsection \ref{Subsec: techlemmas} hold.
Let $a>0$ be given. Assume for a while that typical paths are such that there exists $\epsilon'>0$ such that, $\forall k \in I_N$,
\begin{eqnarray}
&& \max_{k \in I_N} \mathbb{E}_k \left ( \exp{\{a \frac{\max x(t)}{N^{1/3}}\}}\right)<b,\text{ for some $b>0$ },\crcr
&&\nu_{N,o}<N^{1/3-\epsilon'}\text{ and }\nu_{N,e}< \sqrt{\gamma_N}N^{1/3-\epsilon'}.\label{nunty}\end{eqnarray}
The above statement will be proved in the subsequent subsection (Lemma \ref{Lem: max xo } and Lemma \ref{Lem: nun2}) and using exactly the same arguments as in Lemma 7.8 in \cite{PecheFeral}.
We then copy the arguments of Proposition \ref{Prop:z3} and the sequel.
Then it is easy to deduce that typical paths have a non-marked origin, vertices of type $3$ at most 
(and a  number of vertices of type $3$ smaller than $\ln \ln N$) and no edge passed more than twice.
The other paths lead to a negligible contribution. We can also deduce that non-MP-closed vertices of type 2 as well as vertices of type $3$ occur only on the bottom line in typical paths. \\
In particular, let $Z_{3,\infty}$ denote the contribution of paths for which 
$q_1+q_2=0,$ $r_2=0,$ $\sum_{i\geq 4}n_i+\sum_{i\geq 3}p_i=0$, $n_3\leq \ln \ln N$ and no edges read more than twice. 
Then $\mathbb{E}\left [ \text{Tr} M_p^{s_N} \right ]=Z_{3, \infty}(1+o(1))$ and $$Z_{3, \infty}\leq C_6 \sum_{k\in I_N}N \gamma_N^{k-s_N} \sigma^{2s_N}\mathbf{N}(s_N, k)\leq \frac{D_1}{2}v_+^{s_N}.$$
This ensures that the limiting expectation depends only on the variance $\sigma^2$ and has the same behavior as for Wishart ensembles. 
The proof of Theorem \ref{theo: unigammainfini2} is complete, provided we prove the announced Lemmas. $\square$ 

\subsection{Technical Lemmas \label{subsec: techlemaginf}}
We now state the counterpart of Lemma \ref{Lem: max x(t)}.
\bl \label{Lem: max xo } Given $a>0$, there exists $b>0$ such that 
$$\max_{k \in I_N}\mathbb{E}_k \left ( \exp{\{a \frac{\max x(t) (s_N-k)}{N}\}}\right) <b.$$\el
\brem Lemma \ref{Lem: max xo } also yields that
$ \max_{k \in I_N}\mathbb{E}_k \left ( e^{\{a \frac{\max x(t) k}{p}\}}\right)-1 <<1.$
\erem

\paragraph{Proof of Lemma \ref{Lem: max xo }}The proof refers to the proof of Lemma \ref{Lem: max x(t)} in Subsection \ref{Subsec: techlemmas}. Let $l$ be the number of $UU$ steps of a trajectory $X_k,k \in I_N$.
From (\ref{estiml}), one deduces that $l\leq s_N-k\leq A \frac{s_N}{\sqrt {\gamma_N}}.$ Thus by (\ref{estxsos}), we deduce that, if $r\geq a_o \sqrt{As_N \gamma_N^{-1/2}}$, and for any $k \in I_N,$
$\mathbb{P}_k ( \max x(t) =r) \leq \frac{a_3}{\sqrt{s_N-k}}\exp{\{-\frac{a_4r^2}{(s_N-k)}\}}.$
This readily proves Lemma \ref{Lem: max xo }. $\square$
 
\paragraph{}
One next turns to establish the counterpart of Lemma \ref{Lem : nuN}. We denote
$\Gamma_{\nu_{N,o}(X_k)}$ (resp. $\Gamma_{\nu_{N,e}(X_k)}$) the event that the maximal number of times the trajectory $X_k$ comes from above (without falling below) to
some level $x_o$ at even instants (resp. odd instants) is  $\nu_{N,e}$ (resp. $\nu_{N,o}$).
 \bl 
\label{Lem: nun2}There exist positive constants $B_1, B_2, B_3,B_4,$ independent of $N$ and $ p$ such that
\begin{eqnarray}
&&\max_{s_N(1-\frac{A}{\sqrt{\gamma_N}})\leq k \leq s_N(1-\frac{\epsilon}{\sqrt{\gamma_N}})}
\: \mathbb{P}_k(\Gamma_{\nu_{N,o}(X_k)})\leq B_1\frac{s_N^2}{\gamma_N} \exp{\{-B_2 \nu_{N,o}\}}.\label{nuno}\\
&&\max_{s_N(1-\frac{A}{\sqrt{\gamma_N}})\leq k \leq s_N(1-\frac{\epsilon}{\sqrt{\gamma_N}})}
\: \mathbb{P}_k(\Gamma_{\nu_{N,e}(X_k)})\leq B_3\frac{s_N^2}{\gamma_N} \exp{\{-\frac{B_4 \nu_{N,e}}{\sqrt {\gamma_N}}\}}.\label{nune}
\end{eqnarray}
\el

\paragraph{Proof of Lemma \ref{Lem: nun2}}As in the preceding section, we have to estimate the
probability that the trajectory comes to some level $x_o$ $\nu'_N$ times without falling below in some time interval $[s_1,s_2].$ Note that the two steps leading and starting at $s_1$ (resp. $s_2$) are up (resp. down) steps. This is because $[s_1, s_2]$ is a maximal interval. Thus the number of possible choices for $s_1$ and $s_2$ is at most of order $s_N^2/\gamma_N$, as $k\in I_N.$\\
We first prove (\ref{nuno}) and thus assume that $s_1$ is odd. Then the returns to $x_o$ occur at odd instants.
The counterpart of formula (\ref{estimeenu'ninit}) states
\begin{equation}
\mathbb{P}_k\left ( X_k|_{t \in [s'_1, s'_2]} \text{ has only }\eta_N\text{ $UD$ steps }\right ) 
\leq C'_9 \left ( \frac{2A}{\sqrt{\gamma_N}}\right)^{\eta_N}.
\end{equation}
The two steps preceding (and following) $[s'_1,s'_2]$ are either both up steps or both down steps (regardless of the fact that $s'_1, s'_2$ are even or odd). 
The estimate (\ref{estimeenu'n2}) still holds (up to the change $s_N^2\rightarrow s_N^2/\gamma_N$) so that formula (\ref{nuno}) is proved,
copying the proof of Lemma \ref{Lem : nuN}.\\
We now turn to the proof of (\ref{nune}) which is more involved than in Lemma \ref{Lem : nuN}. 
Formula (\ref{estimeenu'ninit}) translates to
$$\mathbb{P}_k \left (X_k \text{ has }\eta_N\text{ UD steps in between } [s'_1,s'_2]\right )
\leq C_9\exp{\{-\epsilon \frac{\eta_N}{\sqrt{\gamma_N}}\}}.$$
Step 1 and Step 2 are then obtained as in Lemma \ref{Lem : nuN} (with $s_N^2\rightarrow s_N^2/\gamma_N$). 
From that, we can deduce that we can consider in Step 3 only the paths for which $A_1\leq Q\leq A_o\nu'_N/\sqrt{\gamma_N}$ for some constants $A_1, A_o>0$. 
We need to refine the estimate for Step $3$. Let then $[s_1,s_2]$ be the interval where $\nu'_N$ returns to some level $x_o$ occur. We call $Y_o$ the trajectory defined by $x(t)-x_o, t\in[s_1, s_2].$ 
We then define $k_o$ to be its number of odd up steps, $Q$ to be its number of returns to $0$ using $DD$ steps, $l$ (resp. $\mu_o$, $\nu_o''$) to be its number of $UU$ steps (resp. of $DU$ steps and of $UD$ steps occuring at some positive level). Assume that $l, Q, \mu_o, k_o$ are given and observe that $k_o=l+\nu'_N-Q+\nu_o''$. 
Let then $\mathbb{P}_{l,Q,k_o,\mu_o}$ denote the conditional probability on the event that $Y_o$ has $k_o$ odd up steps, $\mu_o$ $DU$ steps and $Dy(Y_o)$ has $2l$ steps and $Q$ returns to $0$. Then, one has that 
$$\mathbb{P}_{l,Q,k_o, \mu_o}(\text{ $Y_o$ has $\nu'_N-Q$ $UD$ steps at level $0$ })\leq \frac{C_{\nu'_N}^{Q}C_{(s_2-s_1)/2-\nu'_N}^{\nu_o''}}{C_{(s_2-s_1)/2}^{k_o-l}}.$$
One first shows that it is enough to consider the subpaths $Y_o$ such that $\frac{2k_o}{s_2-s_1}=\frac{k}{s_N}(1+o(1))\in [1-\epsilon_1^{-1}\gamma_N^{-1/2},1-\epsilon_1\gamma_N^{-1/2}]$ for some $\epsilon_1>0$ small enough. This follows from the fact that $k\in I_N$ and arguments already used in Subsection \ref{subsec : clt} (see also (\ref{decroissancel2})). This yields that 
\begin{eqnarray*}&&\mathbb{P}_k\left (\Gamma_{\nu_{N,e}(X_k)}\cap \{\frac{2k_o}{(s_2-s_1)}\leq 1-\epsilon_1^{-1}\gamma_N^{-1/2}\}\right )\crcr
&&\leq \sum_{s_1\leq s_2}\sum_{\frac{2k_o}{(s_2-s_1)}\leq 1-\epsilon_1^{-1}\gamma_N^{-1/2}}\frac{\mathbf{N}(\frac{s_2-s_1}{2},k_o)\mathbf{N}(s_N-\frac{s_2-s_1}{2},k-k_o)}{\mathbf{N}(s_N,k)}\crcr
&&\leq \frac{s_N^2}{\gamma_N} \exp{\{-\eta N^{2/3} \}},\end{eqnarray*}
for some $\eta>0$ provided $\epsilon_1^{-1}>2A$, where $A$ has been fixed in Proposition \ref{Prop: Zinfini0}. The analysis of the case where $2k_o/(s_2-s_1))\geq 1-\epsilon_1\gamma_N^{-1/2}$ is similar. We can assume that $B_4$ is small enough to ensure that $B_4s_N<\sqrt{\gamma_N}\eta N^{2/3}.$ This yields (\ref{nune}) and  ensures that it is enough to consider the case where $\frac{2k_o}{s_2-s_1}=\frac{k}{s_N}(1+o(1))\in [1-\epsilon_1^{-1}\gamma_N^{-1/2},1-\epsilon_1\gamma_N^{-1/2}].$
Fixing $l$ and $k_o$ we set $Q_T:=(1-\frac{k_o-l}{(s_2-s_1)/2})\nu'_N.$ As $l\leq \frac{s_2-s_1}{2}-k_o\leq \frac{\epsilon_1^{-1}}{\sqrt{\gamma_N}}\frac{s_2-s_1}{2},$ one has that $\frac{\epsilon_1}{\sqrt{\gamma_N}}\leq \frac{Q_T}{\nu'_N}\leq \frac{2\epsilon_1^{-1}}{\sqrt{\gamma_N}}.$ 
One can check that there exist a constant $C_o$ independent of $s_1,s_2,k_o$ and $l$ such that, for a given constant $A_2>4$, 
\begin{eqnarray*}&&\!\!\!\!\!\frac{C_{\nu'_N}^{Q}C_{(s_2-s_1)/2-\nu'_N}^{\nu_o''}}{C_{(s_2-s_1)/2}^{k_o-l}}\leq C_o \exp{\{-\frac{(Q-Q_T)^2}{(A_2+1)Q_T}\}}, \text{ if } Q\geq Q_{A_2}:=Q_T(1-A_2),\crcr
&&\!\!\!\!\!\frac{C_{\nu'_N}^{Q}C_{(s_2-s_1)/2-\nu'_N}^{\nu_o''}}{C_{(s_2-s_1)/2}^{k_o-l}}\leq C_o
\exp{\{-\frac{A_2Q_T}{2}\}}\left (\frac{1}{5}\right )^{Q_T(1-A_2)-Q}, \text{ if } Q\leq Q_{A_2}.
\end{eqnarray*}
Thus it is clear that the proportion of paths coming back $\nu'_N$ times from above to some level $x_o$ and for which $Q\leq Q_T(1-\epsilon_1)$
is at most of order $s_N^2/\gamma_N \exp{\{-\frac{\epsilon_1^{3}\nu'_N}{(A_2+1)\sqrt{\gamma_N}} \}}.$
Paths for which $Q\geq Q_T(1-\epsilon_1)\geq \frac{\epsilon_1\nu'_N}{2\sqrt{\gamma_N}}$ are considered as in Step 2. This is enough to ensure (\ref{nune}). $\square$
\brem \label{rem: vargammainf}
The investigation of higher moments follows the same steps as in Subsection \ref{subsec : clt}. In particular, considering $Var(\text{Tr}M_p^{s_N})$, only pairs of correlated paths such that $\mathcal{P}_1\vee \mathcal{P}_2$ has a number of odd up steps of order $2s_N(1-O(\sqrt{\gamma_N^{-1}}))$ are non-negligible.
In (\ref{contrcorrel}), one can also replace the term $(2s_N-1-l)$ with $s_N-(1+l)/2-k_1$. Considering as above the exponential decay of (\ref{majoinfini}), one can also show that
$\sum_{k_1\leq k \wedge 2s_N-1-l}f(k_1)
\leq C'_3(2s_N-1-l-\tilde k_1)^{1/2}f(\tilde k_1)$ where $\tilde k_1= [k \frac{2s_N-1-l}{4s_N-2}](+1).$
Thus (\ref{majosoml}) can be replaced with 
\begin{eqnarray*}
&&\sum_{l=0}^{2s_N-1}C_3(s_N-\frac{l+1}{2}-\tilde k_1)\sqrt{\frac{(2s_N-1-l)}{\sqrt{\gamma_N}}}f(\tilde k_1)\leq C_4 \frac{\sqrt{s_N}\gamma_N^{-3/4}}{(1-\alpha_N)^2},\end{eqnarray*}
where $\alpha_N=k/(2s_N-1)\sim 1-1/\sqrt{\gamma_N}.$ One can readily deduce from the above that the contribution of (\ref{contrcorrelodd}) is negligible. The case where $x(t_e)=0$ yields a negligible contribution, as readily seen from (\ref{nunty}) and Lemma \ref{Lem: nun2}.
The latter is then enough to ensure that $Var(\text{Tr}M_p^{s_N})$ is bounded and only depends on the variance $\sigma^2$ of the entries. The investigation of higher moments is similar.
\erem

%\addcontentsline{toc}{chapter}{Bibliographie}
%\markboth{BIBLIOGRAPHY}{BIBLIOGRAPHY}
%\begin{thebibliography}{nathbib}

\end{document}